
\def\LaTeX{%
  \let\Begin\begin
  \let\End\end
  \let\salta\relax
  \let\finqui\relax
  \let\futuro\relax}

\def\UK{\def\our{our}\let\sz s}
\def\USA{\def\our{or}\let\sz z}

\UK



\LaTeX

\USA


\salta

\documentclass[twoside,12pt]{article}
\setlength{\textheight}{24cm}
\setlength{\textwidth}{16cm}
\setlength{\oddsidemargin}{2mm}
\setlength{\evensidemargin}{2mm}
\setlength{\topmargin}{-15mm}
\parskip2mm


\usepackage[usenames,dvipsnames]{color}
\usepackage{amsmath}
\usepackage{amsthm}
\usepackage{amssymb,bbm}
\usepackage[mathcal]{euscript}

\usepackage{cite}
\usepackage{hyperref}
\usepackage{enumitem}

\usepackage[ulem=normalem,draft]{changes}
%
%

%
 
\definecolor{viola}{rgb}{0.3,0,0.7}
\definecolor{ciclamino}{rgb}{0.5,0,0.5}
\definecolor{blu}{rgb}{0,0,0.7}
\definecolor{rosso}{rgb}{0.85,0,0}

\def\betti #1{{\color{blu}#1}}

\def\pier #1{{\color{rosso}#1}}  
\def\jnew #1{{\color{green}#1}}
\def\pcol #1{{\color{rosso}#1}}

\def\pier #1{{#1}}
\def\betti #1{{#1}}
\def\jnew #1{{#1}}
\def\pcol #1{{#1}}




\bibliographystyle{plain}


%
\newtheorem{theorem}{Theorem}[section]

\newtheorem{definition}[theorem]{Definition}

\newtheorem{lemma}[theorem]{Lemma}

\finqui

\def\Beq{\Begin{equation}}
\def\Eeq{\End{equation}}

\def\Bthm{\Begin{theorem}}
\def\Ethm{\End{theorem}}

\def\Brem{\Begin{remark}\rm}
\def\Erem{\End{remark}}

\def\Bdim{\Begin{proof}}
\def\Edim{\End{proof}}
\def\Bcenter{\Begin{center}}
\def\Ecenter{\End{center}}
\let\non\nonumber




\def\step #1 \par{\medskip\noindent{\bf #1.}\quad}
\def\jstep #1: \par {\vspace{2mm}\noindent\underline{\sc #1 :}\par\nobreak\vspace{1mm}\noindent}

\def\Lip{Lip\-schitz}

\def\rhs{right-hand side}


\def\bh{{\bf h}}
\def\xih{{\xi^{\bf h}}}
\def\etah{{\eta^{\bf h}}}
\def\psih{{\psi^{\bf h}}}


\def\multibold #1{\def\arg{#1}%
  \ifx\arg\pto \let\next\relax
  \else
  \def\next{\expandafter
    \def\csname #1#1#1\endcsname{{\bf #1}}%
    \multibold}%
  \fi \next}

\def\pto{.}

\def\multical #1{\def\arg{#1}%
  \ifx\arg\pto \let\next\relax
  \else
  \def\next{\expandafter
    \def\csname cal#1\endcsname{{\cal #1}}%
    \multical}%
  \fi \next}


\def\multimathop #1 {\def\arg{#1}%
  \ifx\arg\pto \let\next\relax
  \else
  \def\next{\expandafter
    \def\csname #1\endcsname{\mathop{\rm #1}\nolimits}%
    \multimathop}%
  \fi \next}

\multibold
qwertyuiopasdfghjklzxcvbnmQWERTYUIOPASDFGHJKLZXCVBNM.

\multical
QWERTYUIOPASDFGHJKLZXCVBNM.

\multimathop
diag dist div dom mean meas sign supp .


\def\Accorpa #1#2 #3 {\gdef #1{\eqref{#2}--\eqref{#3}}%
  \wlog{}\wlog{\string #1 -> #2 - #3}\wlog{}}


\def\separa{\noalign{\allowbreak}}

\def\<#1>{\mathopen\langle #1\mathclose\rangle}
\def\norma #1{\mathopen \| #1\mathclose \|}

\def\I2 #1{\int_{Q_t}|{#1}|^2}
\def\IT2 #1{\int_{Q_t^T}|{#1}|^2}
\def\IO2 #1{\norma{{#1(t)}}^2}
\def\ov #1{{\overline{#1}}}
\def\next{\\ & \quad}
\def\CP{{\rm (}${\cal CP}${\rm)}}

\def\intQt{\int_{Q_t}}
\def\intQ{\int_Q}
\def\iO{\int_\Omega}

\def\dtt{\partial_{tt}}
\def\dt{\partial_t}
\def\dn{\partial_{\bf n}}
\def\S{{\cal S}}

\def\X{{\cal X}}
\def\Y{{\cal Y}}
\def\Uh{{\cal U}}

\def\checkmmode #1{\relax\ifmmode\hbox{#1}\else{#1}\fi}


\def\bu{{\bf u}}
\def\bh{{\bf h}}

\def\xih{{\xi^{\bf h}}}

\def\muh{{\mu^\bh}}
\def\phih{{\vp^\bh}}
\def\sigmah{{\sigma^\bh}}
\def\wh{{w^
\bh}}
\def\yh{{y^\bh}}
\def\zh{{z^\bh}}

\def\bl{{\boldsymbol \lambda}}

\def\erre{{\mathbb{R}}}
\def\rz{{\mathbb{R}}}

\def\enne{{\mathbb{N}}}

\def\J{{\cal J}}




\def\genspazio #1#2#3#4#5{#1^{#2}(#5,#4;#3)}
\def\spazio #1#2#3{\genspazio {#1}{#2}{#3}T0}

\def\L {\spazio L}

\def\W {\spazio W}

\def\C #1#2{C^{#1}([0,T];#2)}



\def\Lx #1{L^{#1}(\Omega)}
\def\Hx #1{H^{#1}(\Omega)}

\def\Ldue{\Lx 2}

\def\Huno{\Hx 1}
\def\Hdue{\Hx 2}

\def\Liq{{L^\infty(Q)}}
\def\Lio{{L^\infty(\Omega)}}
\def\Lzq{{L^2(Q)}}




\let\vp\varphi

\def\a{\alpha}	
\def\s{\sigma}  
\def\m{\mu}	    
\def\ph{\varphi}	

\def\h{\mathbbm{h}}
\def\G{{\cal G}}

\let\TeXchi\chi                         
\newbox\chibox
\setbox0 \hbox{\mathsurround0pt $\TeXchi$}
\setbox\chibox \hbox{\raise\dp0 \box 0 }
\def\chi{\copy\chibox}



\def\ubar{\overline{\bf u}}
\def\uebar{\overline{u}_1}
\def\uzbar{\overline{u}_2}

\let\hat\widehat

\def\uad{{\cal U}_{\rm ad}}
\def\UR{{\cal U}_{R}}

\def\bvp{\overline\varphi}

\def\bph{{\ov \ph}}
\def\bm{{\ov \m}}   
\def\bs{{\ov \s}}

\usepackage{amsmath}
\DeclareFontFamily{U}{mathc}{}
\DeclareFontShape{U}{mathc}{m}{it}%
{<->s*[1.03] mathc10}{}

\DeclareMathAlphabet{\mathscr}{U}{mathc}{m}{it}
\Begin{document}


%
\title{Optimal control of  a tumor growth model with hyperbolic relaxation of the chemical  potential}

\author{}
\date{}
\maketitle
\Bcenter
\vskip-1.5cm
{\large\sc Pierluigi Colli$^{(1)}$}\\
{\normalsize e-mail: {\tt pierluigi.colli@unipv.it}}\\[0.25cm]
{\large\sc Elisabetta Rocca$^{(1)}$}\\
{\normalsize e-mail: {\tt elisabetta.rocca@unipv.it}}\\[0.25cm]
{\large\sc J\"urgen Sprekels$^{(2)}$}\\
{\normalsize e-mail: {\tt juergen.sprekels@wias-berlin.de}}\\[.5cm]
$^{(1)}$
{{\small Dipartimento di Matematica ``F. Casorati''}\\
{\small Universit\`a di Pavia}\\
{\small via Ferrata 5, I-27100 Pavia, Italy}\\[.3cm] 
$^{(2)}$
{\small Weierstrass Institute for Applied Analysis and Stochastics}\\
{\small Anton-Wilhelm-Amo-Strasse 39, D-10117 Berlin, Germany}\\[10mm]}
\Ecenter
\Begin{abstract}
\noindent 
In this paper, we study the optimal control of a phase field model for a tumor
 growth model of Cahn--Hilliard type in which  the often assumed parabolic relaxation of the chemical potential is replaced  by a hyperbolic one. Both the cases  when the double-well potential 
governing the phase evolution is of either regular or logarithmic type 
are covered by the analysis. We show the Fr\'echet differentiability of
 the associated control-to-state operator in suitable Banach spaces and establish first-order 
 necessary optimality conditions in terms of a variational inequality  involving the adjoint state
 variables. The necessary optimality conditions are then used to derive sparsity results for the
 optimal controls.
 

\End{abstract}

\salta
\pagestyle{myheadings}
\newcommand\testopari{\sc Colli -- Rocca -- Sprekels}
\newcommand\testodispari{\sc {Control of the hyperbolic relaxation in a tumor growth model}}
\markboth{\testopari}{\testodispari}
\finqui
%


\section{Introduction}
\label{INTRO}
\setcounter{equation}{0}

Let  $\a>0,~\tau>0$, and let $\Omega\subset\erre^3$ denote some open and bounded domain having a smooth boundary $\Gamma=\partial\Omega$ with outward normal $\,{\bf n}\,$ and corresponding outward normal derivative $\dn$.  Moreover, we fix some final time $T>0$ and
introduce for every $t\in (0,T]$ the sets $Q_t:=\Omega\times (0,t)$ 
 and $\Sigma_t:=\Gamma\times (0,t)$,
 where we put, for the sake of brevity, $Q:=Q_T$ and $\Sigma:=\Sigma_T$.
We consider in this paper the following optimal control problem: 

\vspace{3mm}\noindent
(${\cal CP}$) \quad Minimize the  cost functional
\begin{align} 
\non	
\J((\mu,\vp,\sigma),{\bf u})
 =&\,  \frac{b_1}2 \intQ |\vp-\widehat \vp_Q|^2
 + \frac{b_2}2 \iO |\vp(T)-\widehat\vp_\Omega|^2
+ \frac{b_3}2 \intQ |{\bf u}|^2 
 \,+\,\kappa\,\G({\bf u})
 \\
 \label{cost} 
=&: \,J((\mu,\vp,\sigma),{\bf u}) + \kappa\G(\bu)
\end{align}
subject to the state system 
\begin{align}
\label{ss1}
&\alpha\dtt\mu+\dt\ph-\Delta\mu=P(\ph)(\sigma+\chi(1-\ph)-\mu) - \h(\ph)
u_1 \quad&&\mbox{in }\,Q\,,\\[1mm]
\label{ss2}
&\tau\dt\vp-\Delta\vp+F^{\,\prime}(\vp)=\mu+\chi\,\sigma \quad&&\mbox{in }\,Q\,,\\[1mm]
\label{ss3}
&\dt\sigma-\Delta\sigma=-\chi\Delta\vp-P(\ph)(\sigma+\chi(1-\ph)-\mu)+u_2\quad&&\mbox{in }\,Q\,,\\[1mm]
\label{ss4}
&\dn \mu=\dn\vp=\dn\sigma=0 \quad&&\mbox{on }\,\Sigma\,,\\[1mm] 
\label{ss5}
&\mu(0)=\mu_0,\quad \dt\mu(0)=\mu_0',\quad\vp(0)=\vp_0,\quad \sigma(0)=\sigma_0 \quad &&\mbox{in }\,\Omega\,,
\end{align}
\Accorpa\Statesys {ss1} {ss5}
and to the control constraint
\Beq
\label{constraint}
{\bf u}=(u_1,u_2)\in\uad\,,
\Eeq
where $\uad$ is some nonempty, closed, bounded and convex subset of the space
$L^\infty(Q)\times L^2(Q)$ that will be specified later.
Moreover, $b_1\ge 0$, $b_2\ge 0$, $b_3>0$  and $\kappa\ge 0$ are prescribed constants,
$\widehat \ph_Q\in L^2(Q)$ and $\widehat\ph_\Omega\in\Huno$ are given target functions, and  
the functional $\,\G:L^2(Q)\times L^2(Q)\to\erre \,$  is nonnegative, continuous and convex.
A prototypical form for $\G$, which enhances the occurrence of {\em sparsity},  is given by
\Beq
\label{formg}
\G({\bf u})= \|{\bf u}\|_{L^1(Q)}=\int_Q (|u_1|+|u_2|)\,.
\Eeq

The  state system \Statesys\ constitutes a simplified and relaxed version of the four-species thermodynamically consistent model for tumor growth
originally proposed by Hawkins-Daruud et al.\ in \cite{HZO} that additionally includes
 chemoctatic terms. 

Let us briefly review the role of the occurring symbols. The primary variables $\ph, \m, \s$ denote the phase field, the associated chemical potential, and the nutrient concentration, respectively.
Moreover, the additional term $\a\dtt\m$ 
is a hyperbolic regularization of equation \eqref{ss1},
while the term $\tau\dt\ph$ is a viscosity contribution to the Cahn--Hilliard equation.

In the above model equations, there are two functions  
acting as distributed controls in the phase and nutrient equations,
respectively. The control variable $u_2$ can model the supply of either a 
medication or some nutrient, while $u_1$, 
which is nonlinearly coupled to the state variable $\varphi$ in the phase equation 
\eqref{ss1} through the term $\,\h(\ph)u_1$, models the application of a cytotoxic drug into 
the system. Here,  $\h$ is a truncation function in order to have the action restricted
to the spatial region in which the tumor cells are located; for instance,
it can be assumed that $\h(-1)=0, \h(1)=1$, and $0<\h(\ph)<1$ if $-1<\ph<1$. 
We refer to  \cite{GLSS, GARL_1, HKNZ, KL} for  possible choices of $\h$. 
The nonlinearity $P$ denotes a proliferation function, and the positive constant $\chi$
represents the chemotactic sensitivity. Finally, the
nonlinear function $F$  is a double-well potential whose derivative represents the 
thermodynamic force driving the evolution of the system. Typical examples, which will be 
included by our analysis, are given by the regular and  logarithmic potentials, which are defined, in this order, by
\begin{align}
\label{regpot}
&F_{\rm reg}(r)=\frac 14 \left(1-r^2\right)^2 \quad\mbox{for $\,r\in\rz$,}\\[1mm]
\label{logpot}
&F_{\rm log}(r)=
\left\{\begin{array}{ll}
(1+r)\,\ln(1+r)+(1-r)\,\ln(1-r)-k_1 r^2\quad&\mbox{for $\,r\in (-1,1)$}\\
2\ln(2)-k_1\quad&\mbox{for $\,r\in\{-1,1\}$\,,}\\
+\infty\quad&\mbox{for $\,r\not\in[-1,1]$}
\end{array}\right.
\end{align}
where $k_1>1$ so that $F_{\rm log}$ is  nonconvex. For technical reasons, 
we do not consider nonsmooth
potentials like the double obstacle potential in this paper.   
The logarithmic potential $F_{\rm log}$ is of particular  relevance in the 
applications, since it is closely connected with the configurational entropy. 
We observe that the thermodynamic force $F^{\,\prime}_{\rm log}(r)$ becomes unbounded as $r\to\pm 1$, 
which forces the phase variable $\vp$ to attain its values in the physically meaningful
range $(-1,1)$. 

As far as well-posedness is concerned, the above model has been investigated in the recent
paper \cite{CRS1}, in which results concerning existence, uniqueness, regularity, and 
continuous dependence have been established (see, in particular, the 
Theorems~\ref{THM:WEAK} and~\ref{THM:REGU}
stated below in Section 2). In this connection, we remark that there 
exists a large body of literature devoted to well-posedness results for different variants 
of the above model. For an account of these contributions, we refer the interested reader to the
introduction and the references given in \cite{CRS1}.  


In this paper, we focus on the optimal control of the state system \Statesys.
The optimal control of tumor growth phase field systems constitutes an important  direction of research, since it can be applied directly to 
devise and implement strategies for the medical treatment of cancer 
patients. Works on boundary and distributed control appeared in~\cite{CGRS3, GARLR,CRW}. Control strategies incorporating chemotaxis, active transport, variable mobilities, and Keller--Segel dynamics were developed in~\cite{ACGW, CSS1, KL,SW,EG, AS, CGSS}. Further contributions included advanced optimality conditions in~\cite{EK,EK_ADV}, as well as a 
refined analysis of treatment-time optimization and related asymptotics in~\cite{S_b,S_a,SigTime,You}. Singular logarithmic and double obstacle potentials 
were addressed in~\cite{S,S_DQ}, while sparse controls and second-order sufficient 
optimality conditions
were studied in~\cite{CSS2,ST,CST}. 
Well-posedness, regularity, and asymptotic behavior for models relevant to control applications and including chemotaxis {were} developed in~\cite{CSS1}. The
abovementioned results collectively provide a rigorous framework for the design and optimization of therapeutic strategies governed by diffuse–interface tumor-growth models. 

Concerning the hyperbolic relaxation of the chemical potential in the viscous Cahn--Hilliard equation (uncoupled from the nutrient and without mass sources), we refer to the recent contributions~\cite{CS,CS-contr}, which inspired the present work. In~\cite{CS}, well-posedness, continuous dependence, and regularity results were established, along with an analysis of the asymptotic behavior as the relaxation parameter $\alpha$ tends to $0$. A related optimal control problem was studied in~\cite{CS-contr}.

The derivation of optimality conditions for the control problem \CP\ makes it necessary to 
establish differentiability properties of the control-to-state operator associated with the state
system \Statesys. The proof of such differentiability properties, in turn, appears to be possible
only if the state variable $\vp$ attains its value in a proper compact subset of the interior of the
domain of the nonlinearity $F^{\,\prime}$ (the interval $(-1,1)$, in the case of the logarithmic potential 
\eqref{logpot}). In other words, the so-called {\em uniform separation condition} must be 
valid. At this point, a significant difference between regular potentials like \eqref{regpot} and
irregular potentials like \eqref{logpot} becomes apparent: indeed, it follows from the Theorems~\ref{THM:WEAK} and~\ref{THM:REGU} below that, under suitable assumptions
on the other data of the system \Statesys, for regular potentials the uniform separation condition is always valid if the control variable $u_1$ belongs to $L^\infty(Q)$, while in
singular cases like \eqref{logpot} this can only be guaranteed under the stronger condition $u_1\in L^\infty(Q)\cap
L^2(0,T;V)$. As a consequence of this fact, the sets $\uad$ of admissible controls and the 
entire analysis of the problem \CP\ are likely to differ in the two cases. We have chosen to
 treat the two cases in separate sections.
    
The paper is organized as follows: in the following Section 2, we formulate the general setting
of our problem and state well-posedness results for the state system \Statesys\ proved in
\cite{CRS1}.  In Section 3, we then consider the case of regular potentials. We show the 
Fr\'e chet differentiability of the control-to-state operator between suitable Banach
spaces and the well-posedness of the associated adjoint state system, and derive first-order
necessary optimality conditions for \CP. In Section 4, the same program is performed for the
singular case, where many parts of the analysis of Section 3 carry over to this situation.  
The final Section 5 then brings some results concerning the sparsity of optimal controls.

Throughout the paper, we make repeated use of H\"older's inequality, of the elementary Young's inequality
\begin{equation}
\label{Young}
a b\,\le \delta |a|^2+\frac 1{4\delta}|b|^2\quad\forall\,a,b\in\erre, \quad\forall\,\delta>0,
\end{equation}
as well as of the continuity of the embeddings $H^1(\Omega)\subset L^p(\Omega)$ for $1\le p\le 6$ and 
$\Hdue\subset C^0(\overline\Omega)$. Notice that the latter embedding is also compact, while this holds true
for the former embeddings only if $p<6$. 


\section{General setting  and the state system}
\label{STATE}
\setcounter{equation}{0}

In this section, we introduce the general setting of our 
problem and state well-posedness results for the state system \eqref{ss1}--\eqref{ss5}
proved in \cite{CRS1}. 
To begin with, for a Banach space $\,X\,$ we denote by $\|\cdot\|_X$
the norm in the space $X$ or in a power thereof, and by $\,X^*\,$ its dual space. 
The only exception to this rule applies to the norms of the
$\,L^p\,$ spaces and of their powers, which we often denote by $\|\cdot\|_p$, for
$\,1\le p\le \infty$. As usual, for Banach spaces $\,X\,$ and $\,Y\,$ that are contained in the same topological vector space, we introduce the linear space
$\,X\cap Y$, which becomes a Banach space when endowed with its natural norm $\,\|u\|_{X \cap Y}:=
\|u\|_X\,+\,\|u\|_Y$, for $\,u\in X\cap Y$.
Moreover, we introduce the spaces
\begin{align}
  & H := \Ldue \,, \quad  
  V := \Huno\,,   \quad
  W := \{v\in\Hdue: \ \dn v=0 \,\mbox{ on $\,\Gamma$}\}.
  \label{defHVW}
\end{align}
Furthermore, we denote by $\,(\,\cdot\,,\,\cdot\,)$, $\,\Vert\,\cdot\,\Vert$, and $\<\cdot,\cdot>$  the standard inner product 
and  norm in $\,H$, as well as the duality pairing between $V$ and its dual $V^*$. We then have the dense and compact embeddings $\,V\subset H\subset  V^*$, with the standard identification  \,\,$\langle v,w \rangle =(v,w)\,\,$ for
 $v\in H$ and $w\in V$. \pier{We also introduce for $t\in[0,T]$ and elements $w\in L^1(0,T;V^*)$ the notation}
\Beq
\label{convolution}
(1\ast w)(t):=\int_0^s w(s)\,ds\,. \
\Eeq

We now provide precise assumptions for the data of the system \Statesys. 
\begin{enumerate}[label={\bf (A\arabic{*})}, ref={\bf (A\arabic{*})}]
\item \label{const:weak}
	$\alpha,\tau $ and $\chi$ are positive constants.
\item \label{F:weak}
	$F=F_1+F_2$ satisfies: $\,F_2\in C^3(\erre)\,$ has a Lipschitz continuous derivative $\,F_2^{\,\prime}$; 
$\,F_1:\erre\to [0,+\infty]\,$ is convex and lower semicontinuous with $\,F_1(0)=0$, and there are constants $\,-\infty\le r_-<0<r_+\le +\infty\,$ such that the restriction of $F_1$ to $(r_-,r_+)$ belongs to $C^3(r_-,r_+)$  and satisfies   
\Beq
\label{singular}
\lim_{r\searrow r_-} F_1^{\,\prime}(r)=-\infty \qquad\mbox{and}\qquad \lim_{r\nearrow r_+} F_1^{\,\prime}(r)=+\infty\,.
\Eeq
\item \label{P:weak}
	$P \in C^2(\erre)$ is nonnegative, bounded, and $P$ and $P'$ are \Lip\ continuous.
\item \label{h:weak}
	$\h \in C^2(\erre)$ is nonnegative, bounded, and $\h$ and $\h'$ are \Lip\ continuous.
\end{enumerate}
Let us note that both the potentials  \eqref{regpot} and \eqref{logpot} fulfill the condition 
${\bf (A2)}$, where in the latter case we have $(r_-,r_+)=(-1,1)$.
Observe also that {\bf (A2)} implies that the subdifferential $\,\partial F_1$ of $F_1$ is a
maximal monotone graph in $\erre\times\erre$. Moreover, by virtue of the growth 
conditions \eqref{singular}, its effective domain $D(\partial F_1 )$
is the open interval $(r_-,r_+)$. Since the restriction of $F_1$ to $(r_-,r_+)$ belongs to
$C^3(r_-,r_+)$ by {\bf (A2)}, we have for every $r\in D(\partial F_1)$ that
$\partial F_1(r)=\{F_1^{\,\prime}(r)\}$.  Also, since $F_1$
attains its minimum value~$0$ at $0$, it turns out that $0\in D(\partial F_1 )$ and $0\in\partial F_1(0)$. Finally, we observe that the assumptions on $F_2$ imply that $F_2$ grows at most quadratically. 

Next, we  introduce our notion of a solution to \Statesys.

\begin{definition}
\label{DEF:WEAK}
A triple $(\m,\vp,\s)$ is called a solution to the state system \Statesys\ if
\begin{align}
&\mu\in H^2(0,T;V^*)\cap W^{1,\infty}(0,T;H)\cap L^\infty(0,T;V),
\label{pier2-0}
\\
&\vp  \in W^{1,\infty}(0,T;H)  \cap H^1(0,T;V) \cap L^\infty(0,T;W)\cap C^0(\overline Q), \label{pier2-1}
\\ 
	&\vp\in (r_-,r_+) \quad\mbox{a.e. in }\,Q, \label{Juve}
	\\
	&\s \in H^1(0,T; H) \cap C^0([0,T];V) \cap \L2 {W}, \label{pier2-2}
	\\ 
	&F_1^{\,\prime}(\vp) \in L^\infty(0,T;H), \label{pier2-3}
\end{align}
and if $(\m,\ph,\s)$ satisfies 
\begin{align}
	 \label{var:1}& \<\alpha\dtt\mu , v > + \iO\dt\ph \, v 
	+ \iO \nabla \mu \cdot \nabla v
	= \iO P(\ph)(\sigma+\chi(1-\ph)-\mu)v
	-\iO \h(\ph)u_1 v
	\nonumber\\
	& \qquad \hbox{for every $v \in V $ and a.e. in $(0,T)$,}\\[2mm]
	\label{var:2}&\tau\dt\vp-\Delta\vp+F^{\,\prime}(\vp)=
	\mu+\chi\,\sigma, \quad \hbox{a.e. in $\,Q$,}
 	\\
 	\label{var:3}
 	& \dt\sigma -\Delta\sigma
 	=-\chi\Delta\ph 
 	-  P(\ph)(\sigma+\chi(1-\ph)-\mu)
 	+ u_2	 \, \,\,\mbox{ a.e. in \,$Q$,}
\end{align}
as well as 
\begin{align}
	\label{var:4}
	\m(0)=\m_0, \quad
\dt\mu(0)=\mu_0',\quad
	\ph(0)=\ph_0, \quad
	\s(0)=\s_0, \quad \hbox{a.e. in $\Omega$}.
\end{align}
\end{definition}
It is worth noting that the homogeneous Neumann boundary 
conditions \eqref{ss4}
are encoded in the conditions~\eqref{pier2-1} and \eqref{pier2-2} for $\ph$ and $\sigma$ (cf.~the definition of the space $W$), and 
in the variational equality \eqref{var:1} for $\mu$. 
Notice also that the initial conditions~\eqref{var:4} are meaningful, 
because  \eqref{pier2-1} and \eqref{pier2-2} imply that $\ph,  \sigma \in
C^0([0,T];V)$, while, owing to \eqref{pier2-0}, it turns out that
$\mu \in   \C1{V^*}\cap \C0H$ and, consequently, $\dt \mu$ is at least weakly continuous from $[0,T]$ to $H$.

We also observe that our notion of solution is a special case of that given in 
\cite[Def.~2.1]{CRS1}. Indeed, there the solutions were quadruples $(\m,\vp,\xi,\s)$ where 
$\xi\in L^\infty(0,T;H)$ satisfied in place of the identity \eqref{var:2} the inclusion condition
$$
\tau\dt\vp-\Delta\vp +\xi+F_2^{\,\prime}(\vp)=\mu+\chi\s, \quad \xi\in\partial F_1(\vp), \quad
\mbox{a.e. in }\,Q.
$$
But, as mentioned above, then \pier{it follows that} $\,\xi=F_1^{\,\prime}(\vp)$ a.e. in $Q$, whence we obtain our equation \eqref{var:2}.

Concerning the well-posedness of the state system, we have the following result.
\begin{theorem}
\label{THM:WEAK}
Suppose that \ref{const:weak}--\ref{h:weak} are fulfilled,  and let the initial data  satisfy
\begin{align}
	\label{weak:initialdata}
	&\mu_0\in V, \quad \mu_0'\in H, \quad
\s_0 \in V,  
	\quad
	\ph_0 \in W,  
\quad \hbox{with} \quad 
r_-< \ph_0(x) <r_+  \quad \mbox{for all }\,x\in\overline\Omega.
\end{align}
Then the state system \Statesys\ has for every  $\,\bu=(u_1, u_2) \in \Liq\times\Lzq$ 
a unique solution  $(\m,\vp,\s)$ in the sense of Definition \ref{DEF:WEAK},
and there exists a constant $K_1>0$, which depends only on 
the norm $\,\|\bu\|_{\Liq\times\Lzq}$  and the data of the system, such that 
\begin{align}
\label{stability}
&\|\mu\|_{H^2(0,T;V^*)\cap W^{1,\infty}(0,T;H)\cap L^\infty(0,T;V) }
\,+ \,\|\vp\|_{W^{1,\infty}(0,T;H)  \cap H^1(0,T;V) \cap L^\infty(0,T;W)\cap C^0(\overline Q)}
\non\\
&\quad +\,\|\xi\|_{L^\infty(0,T;H)}\,+\,\|\s\|_{H^1(0,T; H) \cap C^0([0,T];V) \cap \L2 {W}}\,\le
\,K_1\,.
\end{align}
Moreover, whenever
 $(\m_i,\vp_i,\s_i)$, $i=1,2$, are two solutions to \Statesys\ associated 
with the data $\bu^i=(u_1^i,u_2^i) \in \Liq\times\Lzq$,
$i= 1,2 $, then we have 
\begin{align}
	& \non
	{ \|\m_1-\m_2\|_{\L\infty H}}
	+ \|1\ast(\m_1-\m_2)\|	_{\L \infty V}
		\\
	& \quad \non
	+ \norma{\ph_1-\ph_2}_{\L\infty H \cap \L2 V}
	+ \norma{\s_1-\s_2}_{\L\infty H \cap \L2 V}
	\\
&
\label{cont:dep:weak}
	\leq
	K_2  \Big(
	 \norma{u_1^1-u_1^2}_{\L2 H}
	+  \norma{u_2^1-u_2^2}_{\L2 H}
	\Big)\,,
\end{align} 
with a constant $K_2>0$ which only depends on the data of the system and the norms
$\|\bu^i\|_{\Liq\times\Lzq}$, $i=1,2$.
\end{theorem}
 \Bdim
Observe that, in view of the continuity of the embedding $W\subset C^0(\overline\Omega)$,
the initial datum $\vp_0$ is continuous in $\overline\Omega$.  The last condition in 
\eqref{weak:initialdata} therefore implies  that $\vp_0$ attains its values in a compact subset
of $(r_-,r_+)$. Hence, by virtue of {\bf (A2)}, we have $F_1(\vp_0)\in L^1(\Omega)$ and $F_1^{\,\prime}(\vp_0)\in H$, so that all of the conditions for the application of Theorem~2.2 in 
\cite{CRS1} are met. The result is then a direct consequence of  
\pier{\cite[Thm.~2.2]{CRS1}}.
 \Edim
 
 The next result provides further regularity for more regular initial data and 
 controls $\bu=(u_1,u_2)$. In particular, it yields a {\em uniform separation condition}, which
 will prove to be  fundamental for the control theory in the case of singular potentials like 
 \eqref{logpot} 
 when $-\infty<r_-<r_+<+\infty$. Observe that in the regular case $(r_-,r_+)=(-\infty,
 +\infty)$ such a
 separation condition is automatically satisfied: indeed, the estimate \eqref{stability} yields,
 in particular, that $\|\vp\|_{C^0(\overline Q)}\le K_1$.
\begin{theorem}
\label{THM:REGU}
Assume that \ref{const:weak}--\ref{h:weak} hold true, and let the initial data fulfill \eqref{weak:initialdata} as well as the additional assumptions
\begin{align}
	\label{strong:initialdata}
	\mu_0\in W,\quad \mu_0'\in V, \quad \s_0 \in L^\infty(\Omega). 
\end{align}
In addition, suppose that $\,\bu=(u_1,u_2)\in\Liq\times\Lzq$ satisfies
\begin{align}
	\label{u:strong}
		\bu=(u_1, u_2) \in L^2(0,T;V) \times L^\infty(0,T;H).
\end{align}
Then the solution $(\m,\ph,\s)$ to \Statesys\ in the sense of Definition \ref{DEF:WEAK}
satisfies
\begin{equation}
\label{betti1}
\sigma\in L^\infty (Q), \quad \mu\in H^1(0,T; V)\cap L^\infty(0,T;W), \quad 
F_1^{\,\prime}(\vp) \in L^\infty(Q)\,.
\end{equation}
Moreover, there are constants $K_3>0$ and $r_-<r_* <r^*<r_+$, which depend only on 
 the data of the system and the norm $\,\|\bu\|_{(\Liq\cap L^2(0,T;V))
 \times  L^\infty(0,T;H)}$, such that
\begin{align}
\label{stabu}
&\|\s\|_{L^\infty(Q)}\,+\,\|\m\|_{H^1(0,T;V)\cap L^\infty(0,T;W)}\,+\,\|
F_1^{\,\prime}(\vp)\|_{L^\infty(Q)}\,\le\,K_3\,,
\\[1mm]
\label{separation}
&r_*\le \vp(x.t)\le r^* \quad\mbox{for all }\,(x,t)\in\overline Q\,.
\end{align}
\end{theorem} 
\Bdim
At first, recall that the range of the continuous function $\vp_0$ lies in a compact subset of
$(r_-,r_+)$ so that, by {\bf (A2)}, $ F_1^{\,\prime}(\vp_0)\in L^\infty(\Omega)$.  Therefore, all of the 
conditions for the application of Theorem~3.1 in \cite{CRS1} are met, whence it follows
that the conditions 
\eqref{betti1} and \eqref{stabu} are valid.  
The existence of suitable constants $r_-<r_*<r^*<r_+$ satisfying \eqref{separation} then
follows from the asymptotic growth assumptions made in \eqref{singular}.
\Edim 
 
\section{The case of regular potentials}
\label{REGPOT}
\setcounter{equation}{0}
In this section, we consider the case of regular potentials \pier{defined on the whole real line. Throughout, we \betti{work under the assumptions of Theorem~\ref{THM:WEAK}}. In this setting, 
assumption \jnew{\ref{F:weak}} implies the following regularity property:}
\Beq
\label{glatter}
F_1 \,\mbox{ and }\,F_2\,\mbox{ belong to }\,C^3(\erre).
\Eeq

\noindent  Then, by Theorem~\ref{THM:WEAK}, the control-to-state operator 
\Beq
\label{defS}
\S: \Liq\times L^2(Q) \ni \bu=(u_1,u_2) \mapsto \S(\bu):=(\m,\ph,\s)
\Eeq
is well defined as a mapping into the space 
\begin{align}
\label{defX}
\X\,=\,& \bigl(H^2(0,T;V^*)\cap W^{1,\infty}(0,T;H)\cap L^\infty(0,T;V)\bigr)\non\\
&\times \bigl(W^{1,\infty}(0,T;H)\cap H^1(0,T;V)\cap L^\infty(0,T;W)\cap 
C^0(\overline Q)\bigr)\non\\
&\times \bigl(H^1(0,T;H)\cap \betti{C^0}(0,T;V)\cap L^2(0,T;W)\bigr)\,,
\end{align}
which is a Banach \pier{space} when endowed with its natural norm. Moreover, $\S$ is locally Lipschitz continuous in the sense of \eqref{cont:dep:weak}, \pier{then in a larger space.} 
We now improve this result. 

\Bthm
\label{THM:CD}
Suppose that the assumptions of Theorem~\ref{THM:WEAK} and \eqref{glatter} are satisfied,  and let
$\bu^i=(u_1^i, u_2^i)\in \Liq\times L^2(Q)$ are given with the associated
solutions $\S(\bu^i)=(\mu_i,\ph_i,\s_i)$, for \,$i=1,2$ \betti{with the same initial data}. Then there is a constant $K_4>0$,
which only depends on the data of the system and on the $(\Liq\times L^2(Q))-$norms
of $\bu^1$ and $\bu^2$, such that
\Beq
\label{contdep}
\|\S(\bu^1)-\S( \bu^2)\|_{\cal X}\,\le\,K_4\bigl(\|u_1^1-u_1^2\|_{L^2(0,T;H)}
+\|u_2^1-u_2^2\|_{L^2(0,T;H)}\bigr)\,.
\Eeq
\Ethm 
\Bdim
We introduce the quantities
\begin{align}
&\mu:=\mu_1-\mu_2, \quad \ph:=\ph_1-\ph_2,
\quad \s:=\s_1-\s_2,
\quad u_1:=u_1^1-u_1^2, \quad u_2:=u_2^1-u_2^2\,, \non
\end{align}
which then satisfy the identities
\begin{align}
 \label{diff1}
& \langle \alpha\dtt\mu , v \rangle + \iO \nabla \mu \cdot \nabla v
	= \pier{{}-\iO\dt\ph v +{}} \iO (P(\ph_1)-P(\ph_2))(\sigma_1+\chi(1-\ph_1)-\mu_1)v
\non\\
&	\quad +\iO P(\ph_2)(\s-\chi\ph-\m)v
-\iO (\h(\ph_1)-\h(\ph_2))u_1^1 v
-\iO\h(\ph_2)u_1 v\nonumber\\
	& \qquad \hbox{for every $v \in V $ and a.e. in $(0,T)$,}\\[2mm]
\label{diff2}
&\tau\dt\vp-\Delta\vp \pier{{}= - F^{\,\prime}(\vp_1)  + F^{\,\prime}(\ph_2) +{}} \mu + \chi\,\sigma \quad \hbox{a.e. in $\,Q$,}
 	\\[2mm]
\label{diff3}
 	& \dt\sigma -\Delta\sigma
 	=-\chi\Delta\ph 
 	-  (P(\ph_1)-P(\ph_2))(\sigma_1+\chi(1-\ph_1)-\mu_1)
\non\\
&\hspace*{22mm} - P(\ph_2)(\s-\chi\ph-\mu)
 	+ u_2	 \, \,\,\mbox{ a.e. in \,$Q$,}
\\[2mm]
\label{diff4}
&\mu(0)=0,\quad \dt\m(0)=0,\quad \ph(0)=0,\quad \s(0)=0,\quad
\mbox{ a.e. in $\,\Omega$}. 
\end{align}

Now recall that $\ph_i\in C^0(\overline Q)$, $i=1,2$. Hence, there is some constant $L>0$ such that 
\Beq
\label{differ1}
|P(\ph_1)-P(\ph_2)| + |\h(\ph_1)-\h(\ph_2)| + |F^{\,\prime}(\ph_1)-F^{\,\prime}(\ph_2)| \,\le \,L|\ph| \quad
\mbox{a.e. in \,$Q$}. 
\Eeq
Note that the constant $L$ depends only on the quantity $R:= \max\,\{\|\ph_1\|_{C^0(\overline Q)}\,,\,\|\ph_2\|_{C^0(\overline Q)}\}$, and thus, according to \eqref{stability}, only on the data of the state system \Statesys\ and the norms $\,\|\bu^i\|_{\Liq\times\Lzq}$, $i=1,2$.
In the following, we denote by $C$  such positive constants. 

\pier{We begin by pointing out that the \rhs\ of \eqref{diff2} is at least 
in $\L2 H$
and that its norm can be estimated with the help of \eqref{differ1} and \eqref{cont:dep:weak}. Moreover, the initial datum for $\ph$ in \eqref{diff4} is null 
\jnew{and therefore belongs} to $V$. Then, by applying a well-known parabolic regularity estimate (see, e.g., \cite[Chapter~3]{LioMag}), we find out that 
\begin{align}
	&\norma{\ph}_{H^1(0,T; H) \cap \L\infty V \cap \L2 W}
    \leq C \norma{L|\ph| + |\mu| + \chi |\sigma|}_{\L2 H}\non\\
    &\quad{}
    \leq C \bigl(
	 \norma{u_1}_{\L2 H} +  \norma{u_2}_{\L2 H} \bigr)\,.
    \label{pier1}
\end{align}
\indent Next, the same procedure can be applied to the parabolic equation \eqref{diff3}. We have to verify that the $L^2(0,T;H)$-norm of the \rhs\ of 
\eqref{diff3} is suitably bounded in $L^2(0,T;H)$. In view of  \ref{P:weak}, \eqref{cont:dep:weak}, and \eqref{pier1}, note that only the second term on the \rhs\ of \eqref{diff3} requires a special  
attention. Thanks to the boundedness  of $\s_1,\ph_1,\mu_1$ in 
$L^\infty(0,T;L^4(\Omega))$, the H\"older inequality, and the continuous 
embedding $V\subset \Lx4$, we have that
\begin{align}
&\norma{(P(\ph_1)-P(\ph_2))\bigl(\s_1+\chi(1-\ph_1)-\m_1\bigr)}_{\L2H}^2\non\\
&\le \int_Q L^2|\ph|^2 \bigl(\chi+|\s_1|+\chi|\varphi_1|+|\mu_1|\bigr)^2\non\\
&\le\,C\int_0^T\Bigl(\|\ph(s)\|^2 +\|\ph(s)\|_4^2\bigl(\|\s_1(s)\|_4^2+\|\ph_1(s)\|_4^2
+\|\mu_1(s)\|_4^2\bigr)\Bigr)\,ds\non\\[2mm]
&\le\,C\,\|\ph\|_{L^2(0,T;V)}^2\,\le\,C\bigl(\|u_1\|_{L^2(0,T;H)}^2
+\|u_2\|_{L^2(0,T;H)}^2\bigr)\,. \non
\end{align}
But then it follows from \eqref{pier1} that
\Beq
\label{pier2}
\|\s\|_{H^1(0,T; H) \cap \L\infty V \cap \L2 W}\,\le\,C\bigl(\|u_1\|_{L^2(0,T;H)}
+\|u_2\|_{L^2(0,T;H)}\bigr)\,.
\Eeq 
\indent We now turn our attention to the hyperbolic equation \eqref{diff1}. 
Recalling the theory in \cite[Chapter~3, Sections~8--9]{LioMag}, the solution 
$\mu$ satisfies $\mu \in C^1(0,T;H)\cap C^0(0,T;V)$ and obeys the associated 
energy identity, obtained by formally choosing $v=\partial_t\mu$ in \eqref{diff1} 
and integrating over $(0,t)$. This identity holds for every $t\in[0,T]$.
Taking~\eqref{diff4} into account and applying Young's inequality repeatedly,
we deduce \jnew{that}
\begin{align}
\label{John} 
&\frac{\alpha}{2}\|\partial_t\mu(t)\|^2
+\frac12\|\nabla\mu(t)\|^2
\le \frac12\int_{Q_t}\!\bigl(|\partial_t\varphi|^2+|\partial_t\mu|^2\bigr)
\nonumber\\
&\quad
+L\int_{Q_t} |\ph|\,
\bigl(\chi+|\s_1|+\chi|\varphi_1|+|\mu_1|\bigr)\,|\partial_t\mu|
\nonumber\\
&\quad
+C\int_{Q_t}\!\bigl(|\s|^2+|\varphi|^2
+|\mu|^2+|u_1|^2+|\partial_t\mu|^2\bigr).
\end{align}
Note that, in order to obtain the last term on the \rhs, we have in
particular exploited the bound $\|u_1^1\|_{L^\infty(Q)}$, which is absorbed 
\jnew{in} the constant $C$, together with \jnew{the}  assumptions~\ref{P:weak} and~\ref{h:weak}.
As for the second term on the \rhs\ of \eqref{John}, which we denote
by $I_1$, we apply H\"older's and Young's inequalities, as well as the fact that $\sigma_1$, $\varphi_1$, and $\mu_1$ are bounded in $L^\infty(0,T;L^4(\Omega))$.
It then follows, upon also invoking \eqref{cont:dep:weak}, that
\begin{align}
\label{Paul}
|I_1|\,&\le \, C\int_0^t \|\ph(s)\|_4 \bigl(\pcol{1+{}}\|\sigma_1(s)\|_4+\|\ph_1(s)\|_4+\|\mu_1(s)\|_4\bigr)\|\dt\mu(s)\|\,ds
\non\\
&\le\,C\intQt  |\dt\mu|^2\,+\,C\,\|\ph\|_{L^2(0,T;V)}^2\non\\
&\le\,C\intQt |\dt\mu|^2 \,+\,C\bigl( \|u_1\|^2_{L^2(0,T;H)}+\|u_2\|^2_{L^2(0,T;H)} \bigr)\,.
\end{align}
Combining \eqref{John} and \eqref{Paul} with \eqref{cont:dep:weak} and \eqref{pier1},} we thus obtain that
\begin{align}
\label{Ringo1}
&\frac{\alpha}2 \|\dt\mu(t)\|^2+\frac 12\|\nabla\mu(t)\|^2 \non\\
&\le\,C\intQt |\dt\mu|^2\,+\,
C\bigl( \|u_1\|^2_{L^2(0,T;H)}+\|u_2\|^2_{L^2(0,T;H)} \bigr)\,.
\end{align}
\pier{Then, an application of the Gronwall lemma leads us to the estimate 
\begin{align}
\label{pier3}
\|\mu\|_{W^{1,\infty}(0,T;H)\cap L^\infty(0,T;V)}
\,\le\,C\bigl(\|u_1\|_{L^2(0,T;H)}+\|u_2\|_{L^2(0,T;H)}\bigr)\,.
\end{align}
In addition, from the previous estimates and a comparison in 
\eqref{diff1} it also follows that
\Beq
\label{pier4}
\|\m\|_{H^2(0,T;V^*)}\,\le\,C\bigl(\|u_1\|_{L^2(0,T;H)}
+\|u_2\|_{L^2(0,T;H)}\bigr)\,.
\Eeq
\indent At this point, we may observe that the \rhs\ of \eqref{diff2} is in $H^1(0,T; H) $ and its time derivative is given by
\begin{align*}
& - F''(\ph_1)\dt\ph_1 + F''(\ph_2)\dt\ph_2 + \dt \mu + \chi\, \dt \sigma
\\
&= -(F''(\ph_1)-F''(\ph_2)) \dt\ph_1 - F''(\ph_2)\dt\ph + \dt \mu + \chi\,\dt \sigma,
\end{align*}
which \jnew{belongs to} $\L2 H$, since $|F''(\ph_1)-F''(\ph_2)| \leq C|\ph|$ by \ref{F:weak} and \eqref{stability}, $\dt\ph_1$ is bounded in
$L^2(0,T;V)$, and $F''(\ph_2)$ is obviously bounded. Moreover, the initial value for $\dt\ph $ resulting from \eqref{diff2} and \eqref{diff4} is still zero on account of $\ph_1(0)=\ph_2(0)= \ph_0 $ (cf.~\eqref{var:4}). Then, applying again the parabolic regularity~\cite{LioMag} at the level of $\dt\ph $, it turns out that $\ph \in H^2(0,T; H) \cap \W{1,\infty} V \cap H^1(0,T; W)$ solves 
\begin{align}
\label{diff5}
&\tau\dtt\vp-\Delta (\dt \vp)  \non \\
&= - (F''(\ph_1)-F''(\ph_2)) \dt\ph_1 -  F''(\ph_2)\dt\ph + \dt \mu + \chi\,\dt \sigma \quad \hbox{a.e. in $\,Q$,}
 	\\[2mm]
\label{diff6}
& \dt\ph(0)=0 \quad
\mbox{a.e. in $\,\Omega$}. 
\end{align}
\noindent
Now, multiplying \eqref{diff5} by $\dt\ph$, integration over $Q_t$ and Young's inequality lead to the estimate
\begin{align}
\label{Mick}
&\frac {\tau}2\|\dt\ph(t)\|^2\,+ \,\intQt  |\nabla\dt\ph|^2\non\\
&\leq \,C\intQt |\ph|\, |\dt \ph_1|\, |\dt \ph| + C\intQt |\dt \ph|^2 
+ \intQt |(\dt\m + \chi\dt\s)\dt\ph| \non\\
&\le\, C\int_0^t  \|\ph(s)\|_4\,\|\dt\ph_1(s)\|_4\,\|\dt\ph(s)\|\,ds
+ C\intQt\bigl(|\dt\mu|^2 +  |\dt\sigma|^2 + |\dt\ph|^2 \bigr),
\end{align}
and for the first term on the \rhs\ we have that 
\begin{align}
\label{Jagger}
&C\int_0^t  \|\ph(s)\|_4\,\|\dt\ph_1(s)\|_4\,\|\dt\ph(s)\|\,ds \non\\
&\le\,C\intQt|\dt\ph|^2\,+\,C\|\ph\|_{\L\infty V}^2 \,\|\dt\ph_1\|_{\L2 V}^2 ,
\end{align}
where we used the continuous embedding $V\hookrightarrow L^4(\Omega)$ and the
fact that $\partial_t\varphi_1\in L^2(0,T;V)$.
Collecting \eqref{Mick} and \eqref{Jagger}, and exploiting the previously
derived estimates \eqref{pier1}, \eqref{pier2}, \eqref{pier3}, we conclude that
%
\begin{align}
\label{Midway1}
\|\ph\|_{W^{1,\infty}(0,T;H)\cap H^1(0,T;V)}
\le\,C\bigl(\|u_1\|_{L^2(0,T;H)}+\|u_2\|_{L^2(0,T;H)}\bigr)\,.
\end{align}
Next,} we infer from \eqref{diff2} and \eqref{differ1} that, a.e. in $Q$, 
\begin{align}
\pier{|\Delta\ph|=|-\tau\dt\ph-F^{\,\prime}(\ph_1)+F^{\,\prime}(\ph_2)+\mu  +\chi\s|
\,\le\,C\bigl(|\dt\ph|+|\ph|+|\mu| + |\s|\bigr) }\,.\non
\end{align}
Therefore, we can conclude from \eqref{Midway1} and  standard elliptic 
estimates that
\Beq
\nonumber
\|\ph\|_{L^\infty(0,T;W)} \,\le\,C\bigl(\|u_1\|_{L^2(0,T;H)}
+\|u_2\|_{L^2(0,T;H)}\bigr)\,,
\Eeq
and it follows from \cite[Sect.~8, Cor.~4]{Simon} and the compactness of 
the embedding $W\subset C^0(\overline\Omega)$ that also
\Beq
\nonumber
\|\ph\|_{C^0(\overline Q)}\,\le\,C\bigl(\|u_1\|_{L^2(0,T;H)}
+\|u_2\|_{L^2(0,T;H)}\bigr)\,.
\Eeq
With this, the proof of the assertion is complete.
\Edim
 
 \subsection{The linearized system}

In the following, we study the differentiability properties of the 
control-to-state operator $\S$. We begin our analysis
with the linearization of the system \eqref{ss1}--\eqref{ss5}. To this end,
let $\ubar=(\uebar, \uzbar)\in\Liq\times L^2(Q)$ be fixed and $(\bm,\bph,\bs)
=\S(\ubar)$ denote the associated solution according to Theorem~\ref{THM:WEAK}.
We then consider for given $(h_1,h_2)\in\Lzq\times L^2(Q)$ the 
initial-boundary value problem
\begin{align}
\label{ls1}
&\langle \alpha\dtt\eta,v\rangle +\iO\dt\psi v  + \iO \nabla\eta\cdot\nabla v
\,=\,\iO P(\bph)\bigl(\xi-\chi\psi-\eta\bigr)v\non\\
&\quad +\iO P'(\bph)\bigl(\bs+\chi(1-\bph)-\bm\bigr)\psi v
-\iO \h(\bph)h_1 v- \iO\h'(\bph)\uebar \psi v\non\\[1mm]
&\qquad \mbox{for all $\,v\in V\,$ and a.e. in }\,(0,T),\\[2mm]
\label{ls2}
&\tau\dt\psi - \Delta\psi + F''(\bph)\psi \,=\,\chi\xi+\eta, \quad\mbox{a.e. in }\,Q,
\\[1mm]
\label{ls3}
&\dt\xi-\Delta\xi\,=\,-\chi\Delta\psi -P(\bph)\bigl(\xi-\chi\psi-\eta\bigr)
-P'(\bph)\bigl(\bs+\chi(1-\bph)-\bm\bigr)\psi + h_2, \non\\
&\qquad \mbox{a.e. in }\,Q,\\[1mm]
\label{ls4}
&\dn\psi=\dn\xi=0\quad \mbox{a.e. on }\,\Sigma,\\[1mm]
\label{ls5}
&\eta(0)=\dt\eta(0)=\psi(0)=\xi(0)=0 \quad\mbox{a.e. in }\,\Omega.
\end{align}
Obviously, \eqref{ls1} is just the variational form of the equation 
\begin{align}
\alpha\dtt\eta+\dt\psi-\Delta\eta\,&=\,P(\bph)\bigl(\xi-\chi\psi-\eta\bigr)
+P'(\bph)\bigl(\bs+\chi(1-\bph)-\bm\bigr)\psi \non\\
&\quad - \h(\bph)h_1-\h'(\bph)\uebar\psi\quad \pier{\hbox{in $\,Q$}},\non
\end{align}
together with the boundary condition $\,\dn\eta=0\,$ on $\,\Sigma$. We have the 
following result.
\Bthm
\label{THM:LIN}
Suppose that the assumptions of Theorem~\ref{THM:WEAK} and \eqref{glatter} are fulfilled, and let 
$\ubar=(\uebar,\uzbar)\in\Liq\times L^2(Q)$ be given with associated
$(\bm,\bph,\bs) =\S(\ubar)\in\X$. Then the linearized system \eqref{ls1}--\eqref{ls5}
has for every increment $\,\bh=(h_1,h_2)\in \Lzq\times\Lzq$ a unique solution $\,(\etah,\psih,\xih)
\in\X$, and there exists a constant $K_5>0$, which depends only on the data of 
the system and the $(\Liq\times L^2(Q))-$norm of $\,\ubar$, such that
\Beq
\label{stabulin}
\|(\etah,\psih,\xih)\|_{\X}\,\le\,K_5\bigl(\|h_1\|_{L^2(0,T;H)}+\|h_2\|_{L^2(0,T;H)}\bigr)
\,.
\Eeq 
In other words, the linear mapping $\,\bh\mapsto  (\etah,\psih,\xih)$ is continuous from
$L^2(0,T;H)\times L^2(0,T;H)$ into $\X$. 
\Ethm

\Bdim
The existence result is proved by means of a Faedo--Galerkin approximation using 
the eigenfunctions $\{e_j\}_{j\in\enne}$ of the eigenvalue problem $\,\,-\Delta e_j
=\lambda_j e_j\,$ in $\Omega$, $ \,\dn e_j=0\,$ on $\Gamma$, as basis 
functions. For the sake of brevity, we avoid here writing the finite-dimensional system explicitly and only perform the relevant a priori estimates formally on the 
continuous system \eqref{ls1}--\eqref{ls5}. We note, however, that these formal
estimates are in any case fully justified on the level of the Faedo--Galerkin approximations.
We also observe that $\,\bph\,$ attains its values in a compact subset of $\erre$,
 so that all of the functions $P(\bph)$, $P'(\bph)$, $\h(\bph)$,
$\h'(\bph)$, and $F''(\bph)$, are bounded in $Q$ by a constant that only depends 
on the data of the system and the $\Liq\times\Lzq$--norm of $\ubar$. In the
remainder of this proof, we denote constants having this property by $C$.

Let $t\in (0,T]$ be fixed. We first add $\iO\eta v$ to both sides of \eqref{ls1} and then
choose $\dt\eta$ as test function, a.e. in $(0,T)$. Integrating  with respect to time
over $[0,t]$, invoking {\bf (A3)} and {\bf (A4)}, 
and employing Young's inequality appropriately, we then obtain that 
\begin{align}
\label{lesti1}
&\frac {\alpha}2\|\dt\eta(t)\|^2+\frac 12 \|\eta(t)\|_V^2\non\\
&\le\,C\intQt \bigl(|\eta|^2+(1+\|\uebar\|_{\Liq}^2)|\psi|^2+|\xi|^2 +|h_1|^2+|\dt\eta|^2  \bigr)+\frac {\tau}4 \intQt|\dt\psi|^2 \non\\
&\quad +  C\intQt | \psi |\,|\dt\eta| \,+
 C\int_0^t \bigl(\|\bm(s)\|_4+\|\bph(s)\|_4+\|\bs(s)\|_4\bigr)\|\psi(s)\|_4
\,\|\dt\eta(s)\|\,ds\non\\
&\le\,C\intQt \bigl(|\eta|^2+|\psi|^2+|\xi|^2 +|h_1|^2+|\dt\eta|^2\bigr)
+\frac {\tau} 4\intQt|\dt\psi|^2 +C\int_0^t\|\psi(s)\|_V^2\,ds\,,
\end{align}
where in the last estimate we have also used the continuity of the embedding $V\subset L^4(\Omega)$  and the fact that $\,\bm,\bph,\bs$ are bounded in $L^\infty(0,T;
L^4(\Omega))$.

Next, we add $\psi$ to both sides of \eqref{ls2}, multiply  by $\dt\psi$, and integrate over $Q_t$. Using Young's inequality, we see immediately that
\begin{align}
\label{lesti2}
\tau\intQt|\dt\psi|^2 +\frac 12\|\psi(t)\|_V^2\,\le\,C\intQt\bigl(
|\eta|^2+|\psi|^2+|\xi|^2\bigr)+\frac {\tau}4\intQt|\dt\psi|^2\,.
\end{align}

Finally, we multiply \eqref{ls3} by $\xi$ and integrate  over $Q_t$. Invoking {\bf (A3)},  
employing Young's inequality, and arguing as above in the derivation 
of \eqref{lesti1}, we arrive at the 
estimate  
\begin{align}
\label{lesti3}
&\frac 12\|\xi(t)\|^2+\intQt|\nabla\xi|^2\non\\
&\le\,\chi\intQt \nabla\psi\cdot\nabla\xi
\,+\,C\intQt\bigl(|\eta|^2+|\psi|^2+|h_2|^2\pier{{}+|\xi|^2{}}\bigr)+C\int_0^t \|\psi(s)\|_V^2\,ds\non\\
&\le\,\frac 12\intQt|\nabla\xi|^2+ C\int_0^t\|\psi(s)\|_V^2\,ds
\,+\,C\intQt\bigl(|\eta|^2+|\psi|^2+|h_2|^2\pier{{}+|\xi|^2{}}\bigr)\,.
\end{align}

At this point, \betti{we add the term $\int_{Q_t}|\xi|^2$ to both sides of the inequality \eqref{lesti3}.}  \pier{Next, we sum estimates~\eqref{lesti1}--\eqref{lesti3} 
and rearrange the resulting terms. We can then apply Gronwall's lemma, 
which yields the estimate}
\begin{align}
\label{lesti4}
&\|\eta\|_{W^{1,\infty}(0,T;H)\cap L^\infty(0,T;V)}
\,+\,\|\psi\|_{H^1(0,T;H)\cap L^\infty(0,T;V)} \,+\, \|\xi\|_{L^\infty(0,T;H)\cap L^2(0,T;V)}\non\\
&\le\,C\bigl(\|h_1\|_{L^2(0,T;H)}+\|h_2\|_{L^2(0,T;H)}\bigr)\,.
\end{align}

Having established the estimate \eqref{lesti4}, we can easily conclude. Indeed,
we first observe that the $L^2(0,T;H)-$norm of the expression $\,\chi\xi+\eta-
F''(\bph)\psi-\pier{\tau}\dt\psi\,$ is bounded by an expression of the form $\,\,C(\|h_1\|_{L^2(0,T;H)}+\|h_2\|_{L^2(0,T;H)})$, whence, invoking standard elliptic estimates, we immediately obtain that
\Beq
\label{lesti5}
\|\psi\|_{L^2(0,T;W)}\,\le\,
C(\|h_1\|_{L^2(0,T;H)}+\|h_2\|_{L^2(0,T;H)})\,.
\Eeq
Next, we observe that, thanks to the fact that
$\bm,\bph,\bs$ are bounded in $L^\infty(0,T;L^4(\Omega))$, and owing to the estimate
 \eqref{lesti4},
 it is easily verified that also the 
$L^2(Q)-$norm of the right-hand side of  \eqref{ls3} is bounded 
by $\,\,C(\|h_1\|_{L^2(0,T;H)}+\|h_2\|_{L^2(0,T;H)})$. Therefore, we can conclude from
standard linear parabolic theory that
\Beq
\label{lesti7}
\|\xi\|_{H^1(0,T;H)\cap L^\infty(0,T;V)\cap L^2(0,T;W)}\,\le\,
C(\|h_1\|_{L^2(0,T;H)}+\|h_2\|_{L^2(0,T;H)})\,.
\Eeq
Moreover, a comparison in \eqref{ls1} reveals that 
\Beq
\label{lesti6}
\|\eta\|_{H^2(0,T;V^*)}\,\le\,C(\|h_1\|_{L^2(0,T;H)}+\|h_2\|_{L^2(0,T;H)})\,.
\Eeq

Finally, we differentiate \eqref{ls2} with respect to $t$, obtaining the identity
\Beq
\label{Bond}
\tau\dtt\psi-\Delta\dt\psi\,=\,\chi\dt\xi+\dt\eta-F''(\bph)\dt\psi
-F'''(\bph)\dt\bph\, \psi\,.
\Eeq
The $L^2(Q)-$norm of the \rhs\ is  bounded by 
$\,C(\|h_1\|_{L^2(0,T;H)}+\|h_2\|_{L^2(0,T;H)})$: this is obviuosly true for the first
three summands, while, owing to \eqref{glatter},  the last term can be estimated as follows:
\begin{align}
&\int_Q|F'''(\bph)|^2|\dt\bph|^2\,|\psi|^2\,\le\,C\int_0^T\|\dt\bph(s)\|_4^2
\,\|\psi(s)\|_4^2\,ds\,\le\,C\|\dt\bph\|_{L^2(0,T;V)}^2\,\|\psi\|^2_{L^\infty(0,T;V)}
\non\\
&\le \,C(\|h_1\|_{L^2(0,T;H)}^2+\|h_2\|_{L^2(0,T;H)}^2)\,.\non
\end{align}
Therefore, testing \eqref{Bond} first by $\dt\psi$ and then by $\,-\Delta\psi$,
we easily obtain from Young's inequality and standard elliptic estimates that
\Beq
\|\psi\|_{W^{1,\infty}(0,T;H)\cap H^1(0,T;V)\cap L^\infty(0,T;W)}\,\le\,
C(\|h_1\|_{L^2(0,T;H)}+\|h_2\|_{L^2(0,T;H)})\,.
\Eeq 
Finally, since the embedding $W\subset C^0(\ov\Omega)$ is compact, we can infer
from \cite[Sect.~8, Cor.~4]{Simon} that also
\Beq
\|\psi\|_{C^0(\ov Q)}\,\le\,C(\|h_1\|_{L^2(0,T;H)}+\|h_2\|_{L^2(0,T;H)})\,.
\Eeq

Combining the above estimates, we have thus shown that the system 
\eqref{ls1}--\eqref{ls5} has a solution $(\eta,\psi,\xi)\in\X$ for which the
inequality \eqref{stabulin} is valid.

It is readily seen that the solution is unique. Indeed, if $(\eta_i,\psi_i,\xi_i)\in\X$,
$i=1,2$, are two solutions, then their difference $(\eta,\psi,\xi)=(\eta_1-\eta_2,
\psi_1-\psi_2,\xi_1-\xi_2)$ satisfies the system \eqref{ls1}--\eqref{ls5} with
$h_1=h_2=0$. It then follows from the estimates shown above that $\,\eta=\psi
=\xi=0$, whence the uniqueness follows. With this, the assertion is completely 
proved.
\Edim

\subsection{Existence of optimal controls}
In the remainder of Section 3, we study the optimal control problem \CP\ for
the case of regular potentials. In addition to the assumptions {\bf (A1)}--{\bf (A4)} and  
\eqref{glatter}, we generally assume for the quantities occurring in the cost functional \eqref{cost}:
\begin{enumerate}[label={\bf (A\arabic{*})}, ref={\bf (A\arabic{*})}, start=5]
\item \label{b:coeff}
$b_1\ge 0$, $b_2\ge 0$, $b_3>0$ and $\kappa\ge 0$ are given constants.
\item \label{targets} $\widehat\vp_Q\in \Lzq$ and $\widehat\vp_\Omega\in V$ are prescribed 
target functions.
\item \label{G:term}
$\G:\Lzq\times\Lzq\to\erre$ is nonnegative, continuous and convex.
\end{enumerate}
We now specify the set $\uad$ of admissible controls: we introduce the control space 
\Beq
\label{defUreg}
\Uh:=\Liq\times\Lzq
\Eeq 
and set
\begin{align}
\label{uadreg}
\uad&:=\{{\bf u}=(u_1,u_2)\in\Uh: \underline u_i \le u_i\le\widehat u_i 
\mbox{ \,a.e. in \,$Q$, \,for $i=1,2$} \}\,,
\end{align}
where the threshold functions satisfy  
\begin{align}
\label{thresholdreg}
\underline u_i, \widehat u_i\in\Liq, \,\,\,i=1,2,  \,\mbox{ and  \,$\pier{\underline {u}_i}\le\widehat u_i$,
a.e. in \,$Q$, \,for} \,i=1,2\,.
\end{align} 
At this point, 
we fix once and for all a constant $R>0$ such that
\Beq
\label{defRreg}
\uad\subset\UR:=\{\bu=(u_1,u_2)\in\Uh:\|\bu\|_{\Uh}<R\}\,.
\Eeq 
We then observe the following fact: the constants $K_1,...,K_5$ constructed in the proofs of
the Theorems~\ref{THM:WEAK}, \ref{THM:REGU}, \ref{THM:CD} and~\ref{THM:LIN} can be chosen independently of the controls, as long as 
the latter belong to $\UR$; the constants then depend only on the data of 
the state system and $R$, and no longer on the special control in $\UR$. In particular, it 
follows from \eqref{stability} that the second solution component $\vp$ is uniformly 
bounded in $\overline Q$ by  the same constant $K_1$ provided the associated control $\bu$ 
belongs to $\UR$. Consequently, there exists a constant $K_6>0$, which depends only on $R$ and the data of the state system, such that
\begin{align}
\label{global1}
\max_{i=0,1,2}\left( \|P^{(i)}(\vp)\|_{C^0(\overline Q)} \,+\,\|\h^{(i)}(\vp)\|_{C^0(\overline Q)}\right)
\,+\,\max_{j=1,2}\,\max_{i=0,1,2,3}\,\|F_j^{(i)}(\vp)\|_{C^0(\overline Q\betti{)}}\,\le\,K_6\,,
\end{align} 
whenever $\vp$ is the second solution component corresponding to some $\bu\in\UR$.

We have the following existence result.
\Bthm
\label{THM:EOC}
Suppose that the conditions of Theorem~\ref{THM:WEAK}, \pier{along with}~\ref{b:coeff}--\ref{G:term}, \eqref{glatter} 
\pier{and \eqref{uadreg}--\eqref{defRreg},} are fulfilled. Then the optimal control problem \CP\ has at 
least one solution.
\Ethm
\Bdim
At first, we observe that the cost functional is nonnegative and therefore bounded from below. 
Since $\uad\not =\emptyset$, we can pick a minimizing sequence $\,\{\bu_n\}_{n\in\enne}
\subset\uad$, that is, we have, with $\,(\mu_n,\vp_n,\sigma_n)=\S(\bu_n)$, $n\in\enne$, 
$$
\lim_{n\to\infty} \J(\S(\bu_n),\bu_n) =\inf_{{\bf v}\in\uad} \J(\S({\bf v}) ,{\bf v}) \ge 0.
$$
Since $\uad$ is a closed, bounded, and convex subset of $\Uh$, we may without loss
of generality assume that  $\,\bu_n\to \bu\,$ weakly star in $\Uh$ for some $\bu
\in\uad$. Moreover, since  the bound \eqref{stability} derived in the proof of Theorem~\ref{THM:WEAK} is valid
with a constant $K_1$ that does not depend on the special control in $\UR$, we can claim
that \,$\{(\mu_n,\ph_n,\sigma_n)\}_{n\in\enne}\,$ is bounded in $\X$.  We may therefore assume
that there is a triple $(\mu,\ph,\s)\in\X$ such that 
\begin{align}
\mu_n &\to \mu \quad\mbox{weakly star in }\,H^2(0,T;V^*)\cap \pier{\W{1,\infty}H \cap{}} L^\infty(0,T;V)\,,\non\\
\ph_n &\to\ph \quad\mbox{weakly star in }\,W^{1,\infty}(0,T;H)\cap H^1(0,T;V)\cap
L^\infty(0,T;W)\non\\
&\hspace*{13.5mm} \mbox{and strongly in }\, C^0(\ov Q)\,,\non\\
\sigma_n&\to\s \quad\mbox{weakly star in }\,H^1(0,T;H)\cap L^\infty(0,T;V)\cap L^2(0,T;W)\,.
\non
\end{align}
But then also
$$
P(\ph_n)\to P(\ph),\quad \h(\ph_n)\to\h(\ph), \quad F^{\,\prime}(\ph_n)\to F^{\,\prime}(\ph),\quad
\mbox{all strongly in }\,C^0(\ov Q)\,.
$$
Therefore, taking the limit as $n\to\infty$ in the state system \pier{\eqref{var:1}--\eqref{var:4}},
written for the pairs $((\mu_n,\ph_n,\sigma_n),\bu_n)$, $n\in\enne$, we easily verify that the triple
$(\mu,\ph,\s)$ satisfies the state system for the control $\bu$. By uniqueness,
we must have $(\mu,\ph,\s)=\S(\bu)$, which means that $((\mu,\ph,\s),\bu)$ is 
an admissible pair for \CP.
 
At this point, we notice that the convex and continuous functional $\G$ is weakly
 sequentially lower semicontinuous in $\Lzq\times\Lzq$, which then also applies to the 
  functional
 ${\bf v}\mapsto \J(\S({\bf v}),{\bf v})$. Hence,
 $\,\,
 \J((\S(\bu),\bu)\le \liminf_{n\to\infty}\,\J(\S(\bu_n),\bu_n)\,,
 $   
which means that $\bu\in\uad$ is an optimal control.
\Edim

\subsection{Differentiability of the control-to-state operator}
In this \pier{subsection}, we are going to at prove a differentiability result for the 
control-to-state operator $\S$ in suitable Banach spaces. To this end, we
introduce the space
\begin{align}
\label{defY}
{\cal Y}&:= \bigl(W^{1,\infty}(0,T;H)\cap L^\infty(0,T;V)\bigr)
\times \bigl(H^1(0,T;H)\cap C^0([0,T];V)\cap L^2(0,T;W)\bigr)\non\\
 &\qquad \times\bigl(H^1(0,T;H)\cap C^0([0,T];V)\cap L^2(0,T;W)\bigr)\,,
\end{align}
which is a Banach space when equipped with its natural norm  $\|\,\cdot\,\|_{\cal Y}$.

We have the 
following result.
\Bthm
\label{THM:DIFF}
Suppose that the assumptions of Theorem~\ref{THM:WEAK}, \eqref{glatter}, and \eqref{defRreg} are fulfilled. Then the
control-to-state operator $\S$ is Fr\'echet differentiable in $\UR$ as a mapping from
${\cal U}$ into ${\cal Y}$. Moreover, for every $\ubar=(\uebar,\uzbar)
\in \UR$ the Fr\'echet derivative $D\S(\ubar)\in {\cal L}(
{\cal U}, {\cal Y})$ is given as follows: for every direction $\bh=
(h_1,h_2)\in{\cal U}$, it holds that $D\S(\ubar)[\bh]=(\etah,\psih,\xih)$
is the unique solution to the linearized system \eqref{ls1}--\eqref{ls5}.
\Ethm
\Bdim
Let $\ubar=(\uebar,\uzbar)\in \UR$ be fixed. By Theorem~\ref{THM:LIN}, the linear  
mapping $\,\bh\mapsto  (\etah,\psih,\xih)\,$ is bounded and therefore continuous
from $L^2(Q)\times L^2(Q)$ into $\X$, hence, a fortiori, also as a mapping from
$\Uh$ into the space ${\cal Y}$. 
We have to show that
\Beq
\label{Frechet}
\lim_{\|\bh\|_{\cal U}\to 0}\,\,\frac {\|\S(\ubar+\bh)-\S(\ubar)-(\etah,\psih,\xih)\|_{\cal Y}}
{\|\bh\|_{{\cal U}}} \,=\,0\,.
\Eeq

 Let $\,(\bm,\bph,\bs)=\S(\ubar)$. Then \pier{there is some $\rho>0$ such that $\ubar+\bh\in {\cal U}_R$  for all $\bh\in{\cal U}$
with  $\|\bh\|_{\cal U}< \rho$. In the following, we just consider
such increments $\bh=(h_1,h_2)\in{\cal U}$ with $\|\bh\|_{\cal U}< \rho$, and we denote by $C>0$ constants
that depend only on $R$, $\rho$} and the data of the state system, but not on the special choice of such increments. For any such $\bh$, 
we define the quantities
\begin{align}
&(\muh,\phih,\sigmah):=\S(\ubar+\bh), \quad \yh:=\muh-\bm-\etah, \quad
\zh:=\phih-\bph-\psih, \non\\
&\wh:=\sigmah-\bs-\xih,\non
\end{align}
where $(\etah,\psih,\xih)\in\X$ is the unique solution to the linearized system
\eqref{ls1}--\eqref{ls5} associated with $\bh=(h_1,h_2)$. Clearly, we have 
$(\yh,\zh,\wh)\in\X$, and we see that \eqref{Frechet} certainly holds true if there
is some $C>0$ such that
\Beq
\label{reicht}
\|(\yh,\zh,\wh)\|_{{\cal Y}}\,\le\,C\|\bh\|_{\Uh}^2\,,
\Eeq
which we are going to prove in the following.
To this end, we first observe that all of the constants $C>0$ appearing in the 
\pier{statements} of Theorem~\ref{THM:WEAK} and Theorem~\ref{THM:CD} depend only on the data of the system 
\eqref{ss1}--\eqref{ss5} and the $\Uh-$norm of the controls. Therefore,
there is some constant $C>0$ such that
\begin{align}
\label{Udo1}
&\|(\muh,\phih,\xih)\|_{\X}\le C \quad\mbox{and} \quad
\|(\muh-\bm,\phih-\bph,\sigmah-\bs)\|_{\X} \le C\|\bh\|_{\Uh}\,,\non\\
&\quad\mbox{for all }\,\bh\in\Uh\,\mbox{ with }\,\|\bh\|_{\Uh}\le \pier{\rho}.
\end{align}
In particular, the estimate \eqref{global1} is valid for $\bvp$ and  $\,\phih$, for all 
$\,\bh\in\Uh\,$  with \,$\|\bh\|_{\Uh}\le \pier{\rho}$.

At this point, we subtract from the equations \eqref{ss1}--\eqref{ss5}
satisfied by $(\muh,\phih,\sigmah)$ the sum of the corresponding equations  
for $(\bm,\bph,\bs)$ and the equations \eqref{ls1}--\eqref{ls5} satisfied by
$(\etah,\psih,\xih)$. A little algebra then shows that the quantity
$(\yh,\zh,\wh)$ is a solution to the system
\begin{align}
\label{OSoleMio1}
&\langle \alpha \dtt\yh,v\rangle +\iO \nabla\yh\cdot\nabla v
\,=\,-\iO \dt\zh v + \iO\Lambda_1^\bh v +\iO\Lambda_2^\bh v\non\\[1mm]
&\qquad\mbox{for all $v\in V$ and a.e. $t\in(0,T)$}, \\[2mm]
\label{OSoleMio2}
&\tau\dt\zh-\Delta\zh\,=\,\chi\wh+\yh+
\Lambda_3^\bh \quad\mbox{\,a.e. in $Q$},\\[1mm]
\label{OSoleMio3}
&\dt\wh-\Delta\wh\,=\,-\pier{\chi} \Delta\zh-\Lambda_1^\bh \quad\mbox{\,a.e. in }Q,
\\[1mm]
\label{OSoleMio4}
&\dn \yh=\dn\wh=0 \quad\mbox{\,a.e. on }\Sigma,\\[1mm]
\label{OSoleMio5}
&\yh(0)=\dt\yh(0)=\zh(0)=\wh(0)=0 \quad\mbox{\,a.e. in }\Omega,
\end{align}
with the quantities
\begin{align}
\label{Kurt1}
\Lambda_1^\bh\,&=\,P(\bph)\bigl(\wh-\chi\zh-\yh\bigr)\non\\
&\quad +\bigl(P(\phih)-P(\bph)-P'(\bph)\psih\bigr)\bigl(\bs+\chi(1-\bph)-\bm
\bigr)\non\\
&\quad +\bigl((P(\phih)-P(\bph)\bigr) \bigl((\sigmah-\bs)-\chi(\phih-\bph)
-(\muh-\bm)\bigr)\,,\\[1mm]
\label{Kurt2}
\Lambda_2^\bh\,&=\, - \bigl(\h(\phih)-\h(\bph)-\h'(\bph)\psih\bigr)
\uebar -\bigl(\h(\phih)-\h(\bph)\bigr)h_1\,,\\[1mm]
\label{Kurt3}
\Lambda_3^\bh\,&=\,-\bigl(F^{\,\prime}(\phih)-F^{\,\prime}(\bph)-F''(\bph)\psih\bigr)\,.
\end{align}
By virtue of Taylor's theorem with integral remainder, it holds \pier{that}
\begin{align}
\label{Taylor1}
&P(\phih)-P(\bph)-P'(\bph)\psih \,=\, P'(\bph)\zh + A^\bh
(\phih-\bph)^2,\non\\[1mm]
&\quad{}\h(\phih)-\h(\bph)-\h'(\bph)\psih\,=\,\h'(\bph)\zh + B^\bh
(\phih-\bph)^2, \non\\[1mm]
&\quad\quad{}F^{\,\prime}(\phih)-F^{\,\prime}(\bph)-F''(\bph)\psih \,=\,F''(\bph)\zh +C^\bh (\phih-\bph)^2,
\end{align}
with the remainders
\begin{align}
&A^\bh:= \int_0^1 (1-s)P''(s\phih + (1-s)\bph))\,ds, \non\\
&B^\bh:=\int_0^1 (1-s)\h''(s\phih + (1-s)\bph))\,ds, \non\\
& C^\bh:= \int_0^1 (1-s)F'''(s\phih + (1-s)\bph))\,ds\,. \non
\end{align}
Noting that all  convex combinations of $\phih$ and $\bph$ satisfy the condition 
\eqref{global1}, we see that
\Beq
\label{Udo3}
\|A^\bh\|_{\Liq}  +\|B^\bh\|_{\Liq} +\|C^\bh\|_{\Liq}\,\le\,C\,.
\Eeq 

Now let $t\in (0,T]$ be fixed. We are going to estimate the $L^2(Q_t)-$norms of $\Lambda_i^\bh$, for $i=1,2,3$. As far
as $\Lambda_1^\bh$ is concerned, we estimate the three summands in the
three lines in \eqref{Kurt1}  separately. At first, we obviously have
\Beq
\label{Heino1}
\intQt |P(\bph)|^2 |\wh-\chi\zh-\yh|^2\,\le\,C\intQt\bigl(|\wh|^2+|\yh|^2
+|\zh|^2\bigr) \,.
\Eeq
Next, invoking \eqref{global1}, \eqref{Udo1}, \eqref{Taylor1},  \eqref{Udo3},
the continuity of the embeddings $V\subset L^4(\Omega)$ and $V\subset L^6(\Omega)$, as well as
the fact that $\bm,\bph,\bs$ are bounded in $L^\infty(0,T;V)$, we can argue as follows:
\begin{align}
\label{Heino2}
&\intQt  |P(\phih)-P(\bph)-P'(\bph)\psih|^2\,|\bs+\chi(1-\bph)-\bm|^2\non\\
&\le\,C\intQt \bigl(|\zh|^2+|\phih-\bph|^4\bigr) \bigl(1+|\bm|^2+|\bph|^2+|\bs|^2
\bigr)\non\\
&\le\, C\int_0^t\|\zh(s)\|^2\,ds\,+\,C\int_0^t\|\phih(s)-\bph(s)\|_4^4\,ds
\non\\
&\quad +C\int_0^t \|\zh(s)\|_4^2\bigl(\|\bm(s)\|_4^2+\|\bph(s)\|_4^2+\|\bs(s)\|_4^2
\bigr)\,ds\non\\
&\quad +\,C\int_0^t \|\phih(s)-\bph(s)\|_6^4\bigl(\|\bm(s)\|_6^2+\|\bph(s)\|_6^2+\|\bs(s)\|_6^2
\bigr)\,ds\non\\
&\le\,C\int_0^t\|\zh(s)\|_V^2\,ds\,+\,C\|\bh\|_{\Uh}^4\,.
\end{align}
In addition, using the continuous embedding $V\subset L^4(\Omega)$ and \eqref{Udo1} once more, we find
\begin{align}
\label{Heino3}
&\intQt|P(\phih)-P(\bph)|^2\,|(\sigmah-\bs)-\chi(\phih-\bph)-(\muh-\bm)|^2
\non\\
&\le\,C\int_0^t\Big(\|\phih-\bph\|_4^2\bigl(\|\sigmah-\bs\|_4^2
+\|\phih-\bph\|_4^2+\|\muh-\bm\|_4^2\bigr)  \Big)(s)\,ds\non\\
&\le\,C\|\bh\|_{\Uh}^4\,.
\end{align}
Hence, combining the estimates \eqref{Heino1}--\eqref{Heino3}, we have shown
that 
\Beq
\label{Lambda1h}
\intQt|\Lambda_1^\bh|^2\,\le\,C\intQt\bigl(|\wh|^2+|\yh|^2\bigr)
+\,C\int_0^t\|\zh(s)\|_V^2\,ds \,+\,C\|\bh\|_{\Uh}^4\,.
\Eeq

In addition, we have, using \eqref{Udo1}, \eqref{Kurt2}, \eqref{Taylor1} and \eqref{Udo3}, and arguing similarly
as above,
\begin{align}
\label{Lambda2h}
&\intQt|\Lambda_2^\bh|^2\,\le\,C\intQt
\bigl( |\h(\phih)-\h(\bph)-\h'(\bph)\psih)|^2|\uebar|^2  + |\h(\phih)-\h(\bph)|^2 
|  h_1|^2\bigr)\non\\
&\le \,C\intQt \bigl(|\zh|^2+|\phih-\bph|^4\bigl)+\int_0^t\|(\phih-\bph)(s)\|_\infty
^2\,\|h_1(s)\|^2\,ds\non\\
&\le\,C\intQt|\zh|^2+C\|\bh\|_{\Uh}^4\,.
\end{align} 

Similar reasoning  also yields that
\begin{align}
\label{Lambda3h}
\intQt|\Lambda_3^\bh|^2 \,\le\,C\intQt \bigl(|\zh|^2+|\phih-\bph|^4\bigr)
\,\le\,C\intQt|\zh|^2+ C\|\bh\|_{\Uh}^4\,.
\end{align}

\pier{About the hyperbolic equation \eqref{OSoleMio1}, we recall that the solution $\yh$ is in $C^1(0,T;H)\cap C^0(0,T;V)$ and satisfies the
energy identity (cf.~\cite[Chapter~3, Sections~8--9]{LioMag})
\begin{align*}
\frac \alpha 2\|\dt\yh(t)\|^2 +\frac 12\iO |\nabla \yh(t)|^2
= \intQt  \bigl(- \dt\zh + \Lambda_1^\bh + \Lambda_2^\bh \bigr) \dt\yh \quad\mbox{for all $t\in[0,T]$}.
\end{align*}
This identity is formally obtained by taking $v=\partial_t\yh$ in \eqref{OSoleMio1} 
and integrating over $(0,t)$. Then, we can add the term  $\,\frac 12\|\yh(t)\|^2=\intQt \yh\dt\yh \,$ to both sides and, by virtue of 
Young's inequality, \eqref{Lambda1h}, and \eqref{Lambda2h}, deduce that}
\begin{align}
\label{Heino4}
&\frac \alpha 2\|\dt\yh(t)\|^2 +\frac 12\|\yh(t)\|_V^2 \,\le\,\frac \tau 4\intQt|\dt\zh|^2
\non\\
&\quad{}+C\intQt\bigl(|\wh|^2+|\yh|^2+|\dt\yh|^2\bigr)
\,+\,C\int_0^t\|\zh(s)\|_V^2\,ds\,+\,C\|\bh\|_{\Uh}^4\,.
\end{align}

Next, we add $\zh$ to both sides of \eqref{OSoleMio2}, multiply by $\dt\zh$, and 
integrate over $Q_t$. It then follows from Young's inequality, thanks to \eqref{Lambda3h},
that 
\begin{align}
\label{Heino5}
&\frac 12 \|\zh(t)\|^2_V\,+\,\tau\intQt|\dt\zh|^2\non\\
&\le\,\frac \tau 4\intQt|\dt\zh|^2\,+\,C\intQt \bigl(|\wh|^2+|\yh|^2+|\zh|^2\bigr)
\,+\,C\|\bh\|_{\Uh}^4\,.
\end{align}

Finally, we add $\wh$ to both sides of \eqref{OSoleMio3}, multiply by $\wh$, and 
integrate over $Q_t$. It then follows, invoking Young's inequality and \eqref{Lambda1h} again, that 
\begin{align}
\label{Heino6}
&\frac 12\|\wh(t)\|^2 \,+\intQt|\nabla\wh|^2 \,\le\,\frac 12\intQt|\nabla \wh|^2+\frac {\chi^2}2\intQt|\nabla\zh|^2\non\\
&\quad{}+C\intQt\bigl(|\wh|^2+|\yh|^2\bigr) +C\int_0^t\|\zh(s)\|_V^2\,ds+C\|\bh\|_{\Uh}^4\,.
\end{align}

At this point, we add the inequalities \eqref{Heino4}, \eqref{Heino5} and \eqref{Heino6}.
Rearranging terms, we find that Gronwall's \pier{lemma} can be applied, which then yields the estimate
\begin{align}
\label{OSoleMio6}
&\|\yh\|_{W^{1,\infty}(0,T;H)\cap L^\infty(0,T;V)}+\|\zh\|_{H^1(0,T;H)\cap L^\infty(0,T;V)}
\non\\[1mm]
&\quad{}+\|\wh\|_{L^\infty(0,T;H)\cap L^2(0,T;V)}\le C\,\|\bh\|_{\Uh}^2\,.
\end{align}
But then, in view of \eqref{Lambda3h} and \eqref{OSoleMio6}, also
\begin{align}
\|\Delta\zh\|_{L^2(0,T;H)} \,=\,\|\tau\dt\zh-\chi\wh-\yh-\Lambda_3^\bh\|_{L^2(0,T;H)}
\,\le\,C\|\bh\|_{\Uh}^2\,,\non
\end{align}
and standard elliptic estimates yield that
\Beq
\label{OSoleMio7}
\|\zh\|_{L^2(0,T;W)}\,\le\,C\|\bh\|_{\Uh}^2\,.
\Eeq

Next, we multiply \eqref{OSoleMio3} by $\,-\Delta\wh\,$ and integrate over $Q_t$.
Then, by Young's inequality, \pier{we deduce that} 
\begin{align}
&\frac 12 \|\nabla\wh(t)\|^2\,+\intQt|\Delta\wh|^2
\,\le\,\frac 12\intQt|\Delta\wh|^2\,+\,C\intQt\bigl(|\Delta \zh|^2+|\Lambda_1^\bh|^2\bigr)\,,
\end{align}
and we can infer from the previous estimates and standard elliptic estimates that 
\Beq
\label{OSoleMio8}
\|\wh\|_{L^\infty(0,T;V)\cap L^2(0,T;W)}\, \le\,C\|\bh\|_{\Uh}^2\,.
\Eeq
Comparison in \eqref{OSoleMio3} then yields that also
\Beq
\label{OSoleMio9}
\|\dt\wh\|_{L^2(0,T;H)} \,\le\,C\|\bh\|_{\Uh}^2\,.
\Eeq
Finally, we \pier{recall} the continuity of the embedding 
$$H^1(0,T;H)\cap L^2(0,T;W)\subset C^0([0,T];V),$$
\pier{that, along with the estimates \eqref{OSoleMio6}--\eqref{OSoleMio9}, completes the proof of \eqref{reicht}. Hence, Theorem~\ref{THM:DIFF} is fully proved.}
\Edim

\Brem
\label{REM1}
It is worth mentioning that in the above proof the actual value of the constant $R$ did not
matter (as long as it is large enough to guarantee that $\uad\subset\UR$). We can therefore
claim that the control-to-solution operator $\S$ is Fr\'echet differentiable 
in the sense of Theorem~\ref{THM:DIFF} on the entire control space $\Uh$. 
\Erem

\subsection{First-order necessary optimality conditions}

In this \pier{subsection}, we derive first-order necessary optimality conditions for \CP\ in the case of regular potentials. 

Note that the same optimal control problem has been studied repeatedly for the 
case when the hyperbolic relaxation term $\,\alpha\dtt\mu\,$ is replaced by
a parabolic relaxation term of the form $\,\alpha\dt\mu$. \pier{In particular, references \cite{CSS1, CSS2} address the case without sparsity, i.e., 
for $\kappa = 0$,  while sparsity terms are} included in \cite{CSSJOTA,ST,ST24}, where the latter paper is concerned with
a slightly simplified state system. 

In the following, we will see that at least the first-order
necessary optimality conditions with sparsity term  established in
\cite{CSS1,CSSJOTA,ST} \pier{for the parabolic relaxation 
remain valid, with the appropriate modifications, in the case of hyperbolic 
relaxation. In order to keep the paper at a reasonable length, however, 
we do not address the derivation of second-order sufficient optimality 
conditions as carried out in  \cite{CSS2,ST24}, since this would require 
considerably more involved analytical arguments.}

At this point, we introduce the {\em reduced cost functionals}  $\widetilde\J$ and  $\widetilde J$ by putting
\Beq
\label{redcost}
\widetilde\J(\bu):=\J(\S(\bu),\bu), \quad \widetilde J(\bu):=J(\S(\bu),\bu), \quad 
\mbox{for }\,\bu\in\Uh\,.
\Eeq
Since, thanks to Theorem~\ref{THM:DIFF}, the control-to-state mapping is Fr\'echet differentiable from $\Uh$
into ${\cal Y}$,  the functional $\widetilde J$ is  a Fr\'echet 
differentiable mapping from $\Uh$ into $\erre$. Therefore, the chain rule shows that, for every $\ubar
=(\overline u_1,\overline u_2)\in \Uh$ and $\bh=(h_1,h_2)\in \Uh$, the Fr\'echet 
derivative $D\widetilde J(\ubar)$ satisfies the identity
\begin{align}
\label{DJ}
D\widetilde  J(\ubar)[\bh]\,&=\,b_1\int_Q(\bph-\pier{\widehat\ph_Q})\,\psih\,+\,b_2\iO
(\bph(T)-\pier{\widehat\ph_\Omega})\,\psih(T)\,+\,b_3\int_Q\ubar\cdot\bh,
\end{align}
where $\,(\bm,\bph,\bs)=\S(\ubar)$ and $\pier{(\etah,\psih,\xih)}$  is the solution
to  the linearized system 
\eqref{ls1}--\eqref{ls5}. Notice that Theorem~\ref{THM:LIN} implies that 
the expression on the \rhs\ of \eqref{DJ} is meaningful
also for every $\bh\in L^2(Q)\times L^2(Q)$, and we therefore may extend the Fr\'echet 
derivative $D\widetilde J(\ubar)\in \Uh^*$ to an element of $(L^2(Q)\times L^2(Q))^*$,
which is still denoted by $D\widetilde J(\ubar)$,
by postulating the identity \eqref{DJ} also for general $\bh\in L^2(Q)\times L^2(Q)$.  
In this way, expressions of the form $D\widetilde J(\ubar)[\bh]$ have a well-defined
meaning also for such arguments $\bh$.

We now derive the announced necessary optimality conditions. We recall at this point that a control $\ubar\in\uad$ is called
{\em locally optimal for \CP\  in the sense of $L^p(Q\betti{)}\times L^p(Q)$} for some $p\in[1,\infty]$ if and only if there is some $\varepsilon>0$ such that 
$$
\widetilde\J(\ubar)\le \widetilde\J(\bu) \quad\mbox{for all }\,\bu\in\uad \,\mbox{ such that }
\,\|\bu-\ubar\|_{L^p(Q)\times L^p(Q)}\le\varepsilon\,.
$$
It is easily seen that every locally optimal control in the sense of $L^p(Q)\times L^p(Q)$ for some
$p\in[1,\infty)$ is also locally optimal in the sense of $\Liq\times\Liq$.

Now assume that  $\overline \bu = (\overline u_1,\overline u_2)\in\uad$ 
is a locally optimal control for \CP\ in the sense of $\Liq\times\Liq$.
Then it is easily seen that the variational inequality
\begin{equation}
\label{var1}
{D \widetilde J(\ubar)[\bu - \overline \bu] + \kappa (\G(\bu) - \G(\overline \bu)) \ge 0 \quad \forall\, \bu \in \uad}
\end{equation}
must be satisfied. Indeed, if \,$\bu\in\uad$ is given and $t\in(0,1)$ is sufficiently small, then we have that
$\ubar+t(\bu-\ubar)\in\uad$ and $\|\ubar+t(\bu-\ubar)-\ubar\|_{\Liq\times\Liq}
\le\varepsilon\,.$ Hence, we can infer from the convexity of $\,\G\,$
that
\begin{align*}
0\,&\le\,\widetilde J (\ubar+t(\bu-\ubar))+\kappa \G(\ubar+t(\bu-\ubar))-
\widetilde J(\ubar)-
\kappa \G(\ubar)\\[1mm]
&\le\,\widetilde J(\ubar+t(\bu-\ubar))-\widetilde J(\ubar)+\kappa\,t(\G(\bu)-\G(\ubar)).
\end{align*}

At this point we note that $\S$, being Fr\'echet differentiable from $\Uh$ into $\Y$, is also
Fr\'echet differentiable from the smaller space $\Liq\times\Liq$ into $\Y$.  
Therefore, dividing by $t>0$, and then taking the limit as $t\searrow0$, \eqref{var1} follows. But \eqref{var1} 
implies that $\ubar$ solves the unconstrained convex minimization problem
$$
\min_{\bu \in ( L^2(Q)\times L^2(Q))} \,\bigl( \Phi(\bu)+ \kappa \G(\bu) + I_{\uad}(\bu)\bigr),
$$
where $\Phi(\bu)=D\widetilde J(\overline \bu)[\bu]$ with the extension of the mapping 
$D\widetilde J(\ubar)$ to $\Lzq\times\Lzq$, and where $I_{\uad}$ denotes the indicator function of $\uad$. Hence,
denoting by the symbol $\,\partial\,$ the subdifferential mapping in 
$L^2(Q)\times L^2(Q)$, we have that
$$
\mathbf{0}\in\partial\bigl(\Phi+\kappa \G+I_{\uad}\bigr)(\ubar).
$$
Then we may infer from the well-known rules for subdifferentials of convex functionals that
\[
 \mathbf{0}\in \{D\widetilde J(\ubar)\}+\kappa\partial \G(\ubar)+\partial I_{\uad}(\ubar).
\]
In other words, there
are $\overline\bl\in\partial \G(\ubar)$ and $\widehat\bl\in\partial I_{\uad}(\ubar)$ such that
\begin{equation}\label{supi}
\mathbf{0}=D\widetilde J(\ubar)+\kappa\overline\bl+\widehat\bl.
\end{equation}
But, by the definition of $\partial I_{\uad}(\ubar)$, we have $\,\widehat\bl[\bu-\ubar]\le 0$ for every
$\bu\in \uad$. Hence, thanks to \eqref{supi},
\[
{0\le D\widetilde J(\ubar)[\bu-\ubar]+\kappa\overline\bl[\bu-\ubar]\quad\forall\,\bu\in\uad.}
\] 
Thus, invoking \eqref{DJ}, and identifying  $\overline\bl$  
with the corresponding element of $L^2(Q)\times L^2(Q)$ 
according to the Riesz isomorphism, we have shown the following result.
\begin{lemma}
\label{LEM:COND}
Assume that the \betti{hypotheses} of Theorem~\ref{THM:WEAK}, \ref{b:coeff}--\ref{G:term}, \eqref{glatter}, \eqref{uadreg}
and \eqref{thresholdreg} are fulfilled.
  If \, $\ubar\in\uad$ is a locally optimal control for \CP\ in the sense of 
  $\Liq\times\Liq$, then there is 
   some $\overline \bl \in \partial \G(\ubar)\subset (L^2(Q)\times L^2(Q))$ such that 
\begin{align} \label{varineq1}
&D\widetilde J(\overline \bu)[\bu - \overline \bu] + \kappa \int_Q \overline\bl
\cdot (\bu- \overline \bu)\non\\
&=\,b_1\int_Q(\bph-\pier{\widehat \ph_Q})\,\psih\,+\,b_2\iO
(\bph(T)-\pier{\widehat \ph_\Omega})\,\psih(T)\,
+\int_Q(b_3\ubar+\kappa\ov\bl)\cdot
(\bu-\ubar) \non\\
& \ge 0 \quad \mbox{for all } \,\bu \in \uad\,,
\end{align}
where \, $(\etah,\psih,\xih)$\, is the solution to the 
linearized system \eqref{ls1}--\eqref{ls5} for $\bh=\bu-\ubar$. 
\end{lemma}

Next, we aim at simplifying the expression 
$D\widetilde  J(\overline \bu)[\bu - \overline \bu]$ in \eqref{varineq1}
by introducing an adjoint state.
To this end, we consider the following adjoint system:
\begin{align}
\label{adj1}
&\langle \alpha\dtt p,v\rangle +\iO \nabla p\cdot\nabla v = \iO q v-\iO P(\bph)(p-r)v
\non\\[1mm]
&\qquad\mbox{for all $v\in V$ and a.e. $t\in(0,T)$}\,,\\[2mm]
\label{adj2}
&-\tau\dt q -\dt p -\Delta q =-\chi\Delta r +P'(\bph)\bigl(\bs+\chi(1-\bph)-\bm\bigr)
(p-r)-F''(\bph)q \non\\[1mm]
&\hspace*{39mm} -\chi P(\bph)(p-r)-\h'(\bph)\ov u_1 p+b_1(\bph-\pier{\widehat\ph_Q})
\quad\mbox{a.e. in }\,Q\,,\\[1mm]
\label{adj3}
&-\dt r-\Delta r=\chi q +P(\bph)(p-r)\quad\mbox{a.e. in }\,Q\,,\\[1mm]
\label{adj4}
&\dn q=\dn r=0 \quad\mbox{a.e. on }\,\Sigma\,,\\[1mm]
\label{adj5}
&p(T)=0, \quad \dt p(T)=0, \quad q(T)=\frac{b_2}{\tau}(\bph(T)-\pier{\widehat\ph_\Omega}),
\quad r(T)=0, \quad\mbox{a.e. in }\,\Omega\,.
\end{align}

We have the following well-posedness result. \pier{We emphasize the 
regularity of the solution component $p$ stated in \eqref{regp}, 
which ensures that $p$ is a strong solution of \eqref{adj1}. 
In particular, \eqref{adj1} can be equivalently rewritten as
\begin{align}
\label{adj1bis}
\alpha \partial_{tt} p - \Delta p 
= q - P(\bph)(p - r)
\quad \text{a.e. in } Q,
\end{align}
supplemented with the boundary condition 
$\partial_n p = 0$ a.e. on $\Sigma$.}
\Bthm
\label{THM:ADJ}
Assume that the conditions of Theorem~\ref{THM:WEAK}, \pier{along with}~\ref{b:coeff}--\ref{G:term}, \eqref{glatter}\pier{, \eqref{uadreg}--\eqref{defRreg},} are fulfilled, and let $\ubar\in\UR$  
be given with associated state $(\bm,\bph,\bs)=\S(\ubar)$. 
Then the adjoint system \eqref{adj1}--\eqref{adj5} has a unique solution $(p,q,r)$
with the regularity
\begin{align}
\label{regp}
&p\in \pier{W^{2,\infty}(0,T;H)}\cap W^{1,\infty}(0,T;V)\cap L^\infty(0,T;W)\,,\\
\label{regq}
&q\in H^1(0,T;H)\cap L^\infty(0,T;V)\cap L^2(0,T;W)\,,\\
\label{regr}
&r\in H^1(0,T;H)\cap L^\infty(0,T;V)\cap L^2(0,T;W)\cap \Liq\,.
\end{align}
\Ethm
\Bdim
The existence of a solution with the requested regularity properties can again be shown
by means of a Faedo--Galerkin approximation using the eigenfunctions 
$\{e_j\}_{j\in\enne}$ 
of the elliptic eigenvalue problem $\,\,-\Delta e_j=\lambda_j e_j\,$ in  $\,\Omega$, 
$\,\dn e_j=0\,$ on \,$\Gamma$, as basis functions. Once more, we avoid writing the 
Faedo--Galerkin system explicitly and derive the relevant a priori estimates only
formally by arguing directly on the continuous system \eqref{adj1}--\eqref{adj5}. The
subsequent estimates, while being only formal, are fully justified on the level of the
Faedo--Galerkin approximations. In the following we denote $\,Q^t:=\Omega\times(t,T)$\,
for every $t\in [0,T)$. 
To begin with, let $t\in[0,T)$ be fixed, set 
$$
g_1:=b_1(\bph-\pier{\widehat\ph_Q}), \quad g_2:=b_2(\bph(T)-\pier{\widehat\ph_\Omega}),
$$ 
and observe that $\,g_1\in L^2(0,T;H)$, while $\,g_2\in V$.
  
We test \eqref{adj1} by $\,-\dt p$, integrate  
over $(t,T)$, and add \,\,$\frac 12\|p(t)\|^2 =-\int_{Q^t} p\dt p\,\,$ to both sides of 
the resulting identity. Then, thanks to the terminal conditions for $p$, and owing to Young's inequality, \pier{we have that}
\Beq
\label{adesti1}
\frac \alpha 2\|\dt p(t)\|^2+\frac 12\|p(t)\|_V^2\,\le C\int_{Q^t} \bigl(|p|^2+|q|^2
+|r|^2+|\dt p|^2\bigr)\,.
\Eeq
Likewise, adding $r$ to both sides of \eqref{adj3}, multiplying the resultant by $-\dt r$,
and integrating over $Q^t$, we obtain the estimate
\Beq
\label{adesti2}
\frac 12 \|r(t)\|_V^2+\int_{Q^t}|\dt r|^2\,\le\,
\int_{Q^t}\bigl(|p|^2+|q|^2+|r|^2\bigr)+ \frac 12 \int_{Q^t}
|\dt r|^2\,.
\Eeq
Finally, we multiply \eqref{adj2} by $\,-q\,$ and integrate over $Q^t$, obtaining
\begin{align}
\label{adesti3}
&\frac \tau 2\|q(t)\|^2+\int_{Q^t}|\nabla q|^2\non\\
&\le\,\frac 1{2\tau}\|g_2\|^2 -\chi\int_{Q^t}\nabla q\cdot\nabla r
\,+\,C\int_{Q^t}\bigl(|p|^2+|q|^2+|r|^2+|\dt p|^2 +|g_1|^2\bigr)\non\\
&\quad +\,C\int_{Q^t}\bigl(|\bs|+|\bph|+|\bm|\bigr)\bigl(|p|+|r|\bigr)\,|q|\,.
\end{align}
Now, \pier{we observe that}
\Beq
\label{adesti4}
-\chi\int_{Q^t}\nabla q\cdot \nabla r\,\le\,\frac 12 \int_{Q^t}|\nabla q|^2\,+\,\frac \chi 2
\int_{Q^t}|\nabla r|^2\,\pier{.}
\Eeq
\pier{Besides,} the last term on the \rhs\ of \eqref{adesti3}, which denote by $I$, can be estimated as
follows:
\begin{align}
\label{adesti5}
|I|\,&\le\,C\int_t^T\bigl(\|\bs(s)\|_4+\|\bph(s)\|_4+\|\bm(s)\|_4\bigr) \bigl(\|p(s)\|_4
+\|r(s)\|_4\bigr)\,\|q(s)\|\,ds\non\\
&\le\,\int_{Q^t}|q|^2\,+\,C\int_t^T\bigl(\|p(s)\|_V^2+\|r(s)\|_V^2\bigr)\,ds\,.
\end{align}

At this point, we \pier{apply \eqref{adesti4} and \eqref{adesti5} to the \rhs\ of \eqref{adesti3}, then we sum the result  to  the inequalities \eqref{adesti1} and \eqref{adesti2}}. Rearranging the 
terms, we then find that Gronwall's lemma (taken backward in time)
can be applied, and we conclude the estimate
\begin{align}
\label{adesti6}
&\|p\|_{W^{1,\infty}(0,T;H)\cap L^\infty(0,T;V)}
+\|q\|_{L^\infty(0,T:H)\cap L^2(0,T;V)} +\|r\|_{H^1(0,T;H)\cap L^\infty(0,T;V)}
\non\\
&\le C\bigl(\|g_1\|_{L^2(0,T;H)}+\|g_2\|\bigr)\,.
\end{align}

Next, we observe that, thanks to \eqref{adesti6}, the \rhs\ of \eqref{adj3} is bounded in
$L^\infty(0,T;H)$. Since also $r(T)=0\in\pier{\Lio}$, it follows from classical parabolic 
regularity theory
(see, e.g., \cite[Thm.~7.1, p.~181]{LSU}) that
\Beq
\label{adesti7}
\|r\|_{H^1(0,T;H)\cap L^\infty(0,T;V)\cap L^2(0,T;W)\cap \Liq}\,\le\,C\,.
\Eeq
Similarly, it is now easy to verify that $q$ solves a parabolic problem in which the \rhs\
is already known to be bounded in $L^2(0,T;H)$. Therefore, and since $\pier{g_2}\in V$, it
follows that
\Beq
\label{adesti8}
\|q\|_{H^1(0,T;H)\cap L^\infty(0,T;V)\cap L^2(0,T;W)}\,\le\,C\,.
\Eeq

Finally, we test \eqref{adj1} by $\,\Delta\dt p$ and integrate with respect to time
over $(t,T)$ where $t\in [0,T)$. Integrating by parts, and using the zero terminal conditions for
$p$ and $\dt p$, we obtain from Young's inequality and the above estimates  that
\begin{align}
&\frac \alpha 2 \|\nabla\dt p(t)\|^2+\frac 12\|\Delta p(t)\|^2\,=\,\int_{Q^t}\bigl(q-
P(\bvp)(p-r)\bigr)\Delta\dt p 
\non\\
&=-\iO \bigl(q-P(\bvp)(p-r)\bigr)(t)\Delta p(t) 
\non\\
&\quad{}-\int_{Q^t}\bigl(\dt q-P(\bvp)(\dt p-\dt r)
- P'(\bvp)\dt\bvp(p-r)\bigr)\Delta p
\non\\
&\le \,C\,+\,\frac 14\|\Delta p(t)\|^2\,+\int_{Q^t}|\Delta p|^2
\non\\
&\quad{}+C\int_t^T \|\dt\bvp(s)\|_4 \bigl(\|p(s)\|_4+\|r(s)\|_4 \bigr)\|\Delta p(s)\|\,ds
\non\\ 
&\le \,C\,+\,\frac 14\|\Delta p(t)\|^2\,+\,C\int_{Q^t}|\Delta p|^2\,,  \non
\end{align}
\pier{by noting that $\dt\bvp$ is bounded in $\L2 V$ by \eqref{stability}.  Then, from Gronwall's lemma it follows that}
\Beq
\label{adesti2025}
\|p\|_{W^{1,\infty}(0,T;V)\cap L^\infty(0,T;W)}\,\le\,C\,.
\Eeq
\pier{Consequently, by integrating by parts in \eqref{adj1} and comparing the terms, we arrive at \eqref{adj1bis} with the regularity $\L\infty H$ for $\dtt p$ and the property}
\Beq
\non
\|p\|_{W^{2,\infty}(0,T;\pier{H})}\,\le\,C\,.
\Eeq
\indent With the above estimates, we have shown that the adjoint system 
\eqref{adj1}--\eqref{adj5} has a solution $(p,q,r)$ with the required regularity.
It is easily seen that it is unique: indeed, if $(p_i,q_i,r_i)$, $i=1,2$, are two such 
solutions, then $(p,q,r)=(p_1-p_2,q_1-q_2,r_1-r_2)$ solves the system 
\eqref{adj1}--\eqref{adj5} with the terms $g_1$ and $g_2$ replaced by zero. We then
infer from  \eqref{adesti6} that $p=q=r=0$, whence the uniqueness follows.
This concludes the proof of the assertion.
\Edim

We are now in a position to improve the impracticable result of Lemma~\ref{LEM:COND}.

\Bthm
\label{THM:GOOD-COND}
Suppose that the conditions of Theorem~\ref{THM:WEAK}, \ref{b:coeff}--\ref{G:term}, \eqref{glatter}, \eqref{uadreg}, and
\eqref{thresholdreg} are fulfilled, and let $\ubar\in\uad$ be a locally optimal control for 
\CP\ in the sense of $\Liq\times\Liq$ with the associated state $(\bm,\bph,\bs)=\S(\ubar)$ and the adjoint state $(p,q,r)$. Then there exists some $\ov \bl\in\partial \G(\ubar)$ such 
that
\Beq
\label{VUG}
\int_Q \bigl(\mathbf{d}_0 +\kappa\ov\bl +b_3\ubar\bigr)\cdot(\bu-\ubar)\ge 0
\quad\mbox{for all }\,\bu\in\uad\,,
\Eeq
where $\,\mathbf{d}_0:=(-\h(\bph)p, r)$.
\Ethm

\Bdim
We set $\bh=(h_1,h_2)=(u_1-\uebar,u_2-\uzbar)=\bu-\ubar$ and \pier{let $(\eta,\psi,\xi)$
denote the solution to the linearized system~\eqref{ls1}--\eqref{ls5} corresponding to this increment $\bh$.}
Then we test \eqref{ls1} by $p$, \eqref{ls2} by $q$, and \eqref{ls3}
by $r$, and add the resulting three equations. Integrating by parts with respect to time 
and space, and using the boundary conditions \eqref{ls4} and \eqref{adj4}, the 
initial conditions \eqref{ls5}, as well as the terminal conditions \eqref{adj5}, we obtain:
\begin{align}
\label{Rechnung}
0\,&=\int_0^T\langle \alpha\dtt p(t),\eta(t)\rangle\,dt -\int_Q\psi\dt p +\int_Q\nabla\eta
\cdot\nabla p\non\\
\separa
&\quad+\int_Q p\bigl[-P(\bph)(\xi-\chi\psi-\eta) -P'(\bph)(\bs+\chi(1-\bph)-\bm)\psi
+\h(\bph) h_1+\h'(\bph)\ov u_1\psi\bigr]\non\\
\separa
&\quad+\pier{\tau}\iO\psi(T)q(T)-\pier{\tau} \int_Q \psi\dt q-\int_Q q\Delta\psi+\int_Q q\bigl[F''(\bph)\psi-\chi\xi-\eta\bigr]
\non\\
\separa
&\quad -\int_Q \xi\dt r+\int_Q(-\Delta\xi+\chi\Delta\psi)r\non\\
\separa
&\quad+\int_Q r\bigl[P(\bph)(\xi-\chi\psi-\eta) +P'(\bph)(\bs+\chi(1-\bph)-\bm)\psi -h_2\bigr]
\non\\
\separa
&=\pier{\int_Q \bigl(\h(\bph)\,p\,h_1-rh_2)
+ b_2\iO\psi(T) (\bph(T)-\pier{\widehat\ph_\Omega})}\non\\
\separa
&\quad +\int_0^T\langle \alpha\dtt p(t),\eta(t)\rangle+\int_Q\nabla p\cdot\nabla\eta
+\int_Q\bigl[P(\bph)(p-r)-q\bigr]\eta\non\\
\separa
&\quad +\int_Q\bigl[-\dt r-\Delta r-P(\bph)(p-r)-\chi q\bigr]\xi\non\\
\separa
&\quad +\int_Q\bigl[-\dt p-\pier{\tau}\dt q-\Delta q+\chi \Delta r+
F''(\bph)q+\chi P(\bph)(p-r)\non\\
&\hspace*{20mm} - P'(\bph)(\bs+\chi(1-\bph)-\bm)(p-r)\pier{{}+\h'(\bph)\ov u_1 p}\bigr]\psi\non\\
\separa
&=\,-\int_Q \mathbf{d}_0\cdot\bh +b_1\int_Q \psi(\bph-\pier{\widehat\ph_Q}) +b_2\iO\psi(T)
(\bph(T)-\pier{\widehat\ph_\Omega})\,,
\end{align}
since $(p,q,r)$ solves the adjoint system. Hence, in view of \eqref{varineq1}, the assertion follows.
\Edim

\Brem
\label{REM2}
If the sparsity functional $\G$ has the special form
\Beq
\non
\G(\bu)=\G_1(u_1)+\G_2(u_2)=\G_1(I_1(\bu))+\G_2(I_2(\bu)) \quad\mbox{for }\,\bu=(u_1,u_2)\in L^2(Q)\times L^2(Q),
\Eeq
with nonnegative, convex and continuous functionals $\G_i:L^2(Q)\to\erre$
and the linear projection mappings $I_i(\bu)=u_i$, for $i=1,2$,
then it follows from the sum and chain rules for subdifferentials 
(cf.,~\pier{e.g.,}~\cite[Sect.~4.2.2, Thm.~1 and Thm.~2]{IT}) that 
$$
\partial \G(\ubar)=\{\bl=(\lambda_1,\lambda_2):\lambda_i\in\partial \G_i(\ov u_i), \,\,\,i=1,2\}.
$$
In this case, we have $\ov\bl=(\ov\lambda_1,\ov\lambda_2)$ with $\ov\lambda_i\in
\partial \G_i(\ov u_i)$ for $i=1,2$, and 
\eqref{VUG} is equivalent to two independent variational inequalities that
have to \pier{be} satisfied simultaneously, namely
\begin{align}
\label{VUG1}
&\int_Q\bigl(-\h(\bph)p+\kappa\ov\lambda_1+b_3\ov u_1\bigr)\bigl(u_1-\uebar\bigr)
\ge 0 \quad\mbox{if }\,u_1\in\Liq \,\mbox{ and }\,\underline u_1\le u_1\le\hat u_1,\\
\label{VUG2}
&\int_Q\bigl(r+\kappa\ov\lambda_2+b_3\uzbar\bigr)\bigl(u_2-\uzbar\bigr)\ge 0
\quad\mbox{if }\,u_2\in\Liq \,\mbox{ and }\,\underline u_2\le u_2\le\hat u_2.
\end{align}
Then a standard argument leads to the pointwise projection formulas
\begin{align}
\label{Pro1}
&\uebar\,=\,\max\,\Big\{\underline u_1,\,\min\,\Big\{\hat u_1,-b_3^{-1}(-\h(\bph)p
+\kappa\ov\lambda_1)\Big\}\Big\} \quad\mbox{a.e. in }\,Q\,,\\
\label{Pro2}
&\uzbar\,=\,\max\,\Big\{\underline u_2,\,\min\,\Big\{\hat u_2, -b_3^{-1}(r+\kappa
\ov\lambda_2)\Big\}\Big\} \quad\mbox{ a.e. in }\,Q\,.
\end{align}
\Erem

\section{The case of singular potentials}
\label{SINGPOT}
\setcounter{equation}{0}

We now turn our interest to nonregular potentials. To this end, we assume the condition:
\Beq
\label{unglatt}
\mbox{\pier{it holds that} }\, -\infty<r_-<0<r_+<+\infty. 
\Eeq  
In this case, the conditions of Theorem~\ref{THM:WEAK} do not suffice to carry out an analysis as in
the previous section. Indeed, although the boundedness condition \eqref{stability} is
valid also in this situation, it fails to guarantee the uniform boundedness
of $\,\|F_1^{\,\prime}(\vp)\|_{\Liq}$, since it cannot be excluded that 
\pier{the values of the phase variable $\vp$
get arbitrarily close} to the critical values $r_-,r_+$ at which $|F_1^{\,\prime}|$ blows up. 
Therefore, the estimate \eqref{global1}, which was
an indispensable prerequisite for the analysis of Section 3 to work, is not at our disposal.
In this connection, notice that also the uniform Lipschitz inequality \eqref{differ1}, which was fundamental
in the proof of Theorem~\ref{THM:CD}, is a consequence of \eqref{global1}.
To overcome this inherent difficulty, we assume for the remainder of this section that the 
stronger conditions of Theorem~\ref{THM:REGU} are fulfilled, since then the uniform separation 
\eqref{separation} is satisfied (which implies the validity of \eqref{global1}), provided that the 
$(\Liq\cap L^2(0,T;V))\times L^\infty(0,T;H))$--norms of the admissible controls are uniformly bounded. \jnew{It therefore follows that under the assumptions of 
Theorem~\ref{THM:REGU} the Fr\'echet differentiability result of Theorem~\ref{THM:DIFF} is also valid in
the singular case when only \eqref{unglatt} is satisfied.}

It is not advantageous  to \jnew{postulate} that a control variable 
($u_1$, in this case) belongs to a bounded subset of $L^2(0,T;V)$, since it imposes a severe restriction on the available control spaces in practical situations. In particular, the usual
techniques to derive pointwise projection formulas from variational inequalities (cf.~Remark~\ref{REM2}) fail. We therefore restrict ourselves to the following special cases.

\vspace{2mm} \noindent
\underline{\em Scenario 1:}  \quad We consider $u_1\in \Liq\cap L^2(0,T;V)$ as a fixed
given datum and apply the control action only to the second variable $u_2$. 
For convenience, we write $\ov u_1$ in place of $u_1$, in the following.
This 
becomes formally a
special case of the one discussed in Section 3 if we make the 
following choice (cf.~\eqref{thresholdreg})
in the definition \eqref{uadreg} of $\uad$:
\Beq
\label{thresholdsc1}
\underline u_1=\ov u_1=\widehat u_1 \pier{{} \in \Liq\cap L^2(0,T;V)}\,.
\Eeq
Clearly, then \CP\ becomes in reality a control problem only in the variable $u_2$,
but we can take advantage of the machinery developed in the previous section. 
Indeed, the uniform  separation condition \eqref{separation} is valid for all second solution components associated with controls $\bu\in\UR$, and the 
fundamental estimate \eqref{global1} holds true. From this point on, it is easily seen
that the whole analysis carried out in  Section 3 can be repeated accordingly. Since in this case we actually have to deal only with the control variable $u_2$, the analysis even simplifies,
 yielding the existence of optimal controls and a first-order 
necessary optimality condition that resembles \eqref{VUG} in Theorem~\ref{THM:GOOD-COND}. 

Notice, in particular, that the linearized system \eqref{ls1}--\eqref{ls5} changes slightly:  in fact,
in the second line of Eq.~\eqref{ls1}, the term $\,\,-\iO \h(\ov \ph)h_1 v\,\,$ does not occur. 
Since the form of the adjoint state system \eqref{adj1}--\eqref{adj5} remains unchanged 
(recall that we write $\ov u_1 $ in place of $u_1$), this implies 
that also in the calculation carried out in \eqref{Rechnung} the corresponding term 
$\,\,\int_Q \h(\ov \ph)\,p\,h_1\,$\, does not occur, which clearly corresponds to the
fact that the first control component is fixed. Accordingly, the first component of
the vector function ${\bf d}_0$ introduced in Theorem~\ref{THM:GOOD-COND} vanishes, and only the part of the variational inequality \eqref{VUG}, which acts on the second control component, yields some 
useful information. To write this condition more
explicitly,  we assume that the sparsity functional $\G$ is of the special form 
\Beq
\label{defGsc1}
\G(\bu)=\G((u_1,u_2))=\G_2(u_2) \quad\mbox{for every }\,\bu=(u_1,u_2)\in \Lzq\times\Lzq\,,
\Eeq
where $\G_2:L^2(Q)\to\erre$ is nonnegative, continuous and convex. 
We also introduce the set
\Beq
\label{uadsc1}
\uad^1:=\{ u\in\Liq: \underline u_2\le u\le \widehat u_2\,\mbox{ a.e. in }\,Q\}.
\Eeq
A locally optimal control for \CP\ in the sense of $\Liq$ is then any  
$u\in\uad^1$ such that there is some $\varepsilon >0$ such that,  with the fixed 
given function $\ov u_1$,
\Beq
\J(\S((\ov u_1,u),(\ov u_1,u))\betti{)}\le \J(\S((\ov u_1,v),(\ov u_1,v))\betti{)} \quad\mbox{for all }\,v\in\uad^1\,
\mbox{ with }\,\|v-u\|_\Liq\le \varepsilon\,.\non
\Eeq 

We then conclude from Theorem~\ref{THM:GOOD-COND} the following optimality condition. 
\Bthm
\label{OPT1}
Suppose that the conditions of Theorem~\ref{THM:REGU} are satisfied with the exception of
\eqref{u:strong}. In addition, let \pcol{\ref{b:coeff}}--\ref{G:term},  \eqref{thresholdreg},
and \eqref{unglatt}--\eqref{uadsc1} be fulfilled. Moreover,
let \,$\overline u_2\in\uad^1$, with associated state
$(\overline \mu,\overline \ph , \overline\s)=\S((\ov u_1,\overline u_2))$ and adjoint state 
$(p,q,r)$, be a locally optimal control for \CP\ in the sense of $\Liq$. 
Then there is some $\ov \lambda_2\in\partial\G_2(\uzbar)$ such that the 
variational inequality \eqref{VUG2} and the pointwise projection formula
\eqref{Pro2} are satisfied. 
\Ethm  

\vspace{2mm}\noindent
\underline{\em Scenario 2:}
\quad Assume that the assumptions Theorem~\ref{THM:REGU} are satisfied, except the condition
\eqref{u:strong}.  We then fix an element $\widehat z\in \pier{\Lio}\cap V$ and 
consider controls $u_1$ of the product type
\Beq 
\label{Sce2.1}
u_1(x,t)= \widehat z(x)u(t) \quad\mbox{for a.e. }\,(x,t)\in Q\,,
\Eeq
that is, the \jnew{real control variable $u$ depends only on time.}
Accordingly, we  consider the set of admissible controls
\begin{align}
\label{uadsc2}
&\uad^2:=\bigl\{(u,u_2)\in L^\infty(0,T)\times\Liq: \underline u_1\le u\le \widehat u_1 \,
\mbox { a.e. in $\,(0,T)$} 
\non\\
& \qquad \qquad \mbox{and }\,\underline u_2\le u_2\le \widehat u_2 \,\mbox{ a.e. in }\,Q  \bigr\}\,,
\end{align}
\pier{where
\begin{align}
\label{thresholdsc2}
&\underline u_1 , \widehat u_1 \in L^\infty(0,T) \ \hbox{ and } \ 
\underline u_1\le\widehat u_1\,\mbox{ a.e. in $\,(0,T)\,$}, \non \\
&\quad 
\underline u_2, \widehat u_2 \in \Liq \ \hbox{ and } \ 
\underline u_2\le\widehat u_2\,\mbox{ a.e. in \,}Q\,.
\end{align}
We} chose some $R>0$ such that 
\begin{align}
\label{defRsc2}
&\uad^2\subset\UR:=\bigl\{(u,u_2)\in L^\infty(0,T)\times\Liq: 
\|(u,u_1)\|_{L^\infty(0,T)\times\Liq} <R \bigr\}\,.
\end{align}
Then, for every $(u,\jnew{u_2})\in\UR$, the function $(\widehat z\,u,
\jnew{u_2})$ belongs to the open
ball of radius $\,\max \bigl\{1, \|\widehat z\|_{\Liq\cap L^2(0,T;V)}\bigr\}R$ \,in the space
$(\Liq\cap L^2(0,T;V))\times\Liq$.

Therefore, the uniform  separation condition \eqref{separation} is valid for all second solution components associated with controls $(u,u_2)\in\UR$, and thus also \eqref{global1}, 
so that the analysis
 carried out in  Section 3 can be repeated analogously. 
 
 \jnew{In this situation, we consider the cost functional in the form}\pcol{%
 \[ \J_2((\mu,\vp,\s),(u,u_2)= J_2((\mu,\vp,\s),(u,u_2)\pier{)}+\kappa \G((u,u_2))\,,\]
with the differentiable part}
 \begin{align}
 \label{costdue}
 J_2((\mu,\vp,\s),(u,u_2))=&\ \frac{b_1}2\int_Q|\vp-\widehat\vp_Q|^2+\frac{b_2}2\iO
 |\vp(T)-\widehat \vp_\Omega|^2 \non \\
 &{}+\frac{b_3}2
 \jnew{\biggl(\int_0^T|u(t)|^2\,dt\,+\int_Q|u_2|^2\biggr)},
 \end{align}
 \jnew{and where we assume the sparsity term  in the additive form
 \begin{align}
 \label{lumpi1}
 {\cal G}((u,u_2))={\cal G}_1(u)+{\cal G}_2(u_2)  \quad\mbox{for 
 \,$(u,u_2)\in L^2(0,T)\times L^2(Q)$,}
 \end{align}
 with nonnegative, continuous and convex functionals ${\cal G}_1:L^2(0,T)\to\erre$ and
 ${\cal G}_2:L^2(Q)\to\erre$. According to Remark~\ref{REM2}, it  then follows that
 \begin{align}
 \partial {\cal G}((u,\pcol{u_2}))=\{(\lambda_1,\lambda_2):\lambda_1\in\partial {\cal G}_1(u)
 \,\,\,\mbox{and}\,\,\,\lambda_2\in \partial {\cal G}_2(u_2)\}. \nonumber
 \end{align}
 Notice also that in \eqref{costdue} we have
 }
 $(\mu,\vp,\s)=\S({\cal H}[u],u_2)$ with the linear operator
 \Beq
\label{Hsc2}
{\cal H}:L^2(0,T)\ni u\mapsto \widehat z\,u=:{\cal H}[u] \in \Lzq\,.
\Eeq 
Clearly ${\cal H}$ is continuous both from $L^2(0,T)$ into $\Lzq$ and from $L^\infty(0,T)$ into $\Liq\cap L^2(0,T;V)$  and, as a linear operator, also Fr\'echet differentiable between these spaces with the Fr\'echet derivative
$D{\cal H }={\cal H}$. It \jnew{therefore follows from the} chain rule that the 
control-to-state operator, which in this situation is given by the operator ${\cal T}$,
where 
\Beq
\label{defTsc2}
{\cal T}((u,u_2)):=\S({\cal H}[u],u_2)=\S(\widehat z\,u,u_2) \quad\mbox{for }\,\,(u,u_2)\in L^\infty(0,T)\times\Liq\,,
\Eeq  
is Fr\'echet differentiable from $L^ \infty(0,T)\times\Liq$ into \jnew{${\cal Y}$}. For fixed $(\ov u,\uzbar)
\in L^\infty(0,T)\times \Liq$, its derivative $D{\cal T}(\ov u,\uzbar)$ is given as follows:
for all directions ${\bf h}=(h,h_2)\in L^\infty(0,T)\times \Liq$, the value $D{\cal T}(\ov u,\ov u_2)
[(h,h_2)]$ is given by the unique solution $(\etah,\psih,\xih)$ to the linearized system 
\eqref{ls1}--\eqref{ls5}, in which in this case the sum of the last two summands in the second line of Eq.~\eqref{ls1}  has to be replaced by the expression 
$$ 
-\iO \h(\bvp)\,\widehat z\,h\,v -\iO \h'(\bvp)\,\widehat z\,\ov u\,\psi\,v\,.
$$
Consequently,  the penultimate term in the second line of Eq.~\eqref{adj2} in the adjoint system 
must be replaced by   $\,\,-\h'(\bvp)\,\widehat z\,\ov u\,p$. 
Accordingly, the sum of the
corresponding last two summands in the second line of the calculation \eqref{Rechnung},
which now is performed for ${\bf h}=(h,h_2)=(u-\ov u,u_2-\uzbar)$,
becomes 
$$  
\int_Q \bigl(\h(\bvp)\,\widehat z\,h \,+\,\h'(\bvp)\,\widehat z\,\ov u\,\psi\bigr)p\,.
$$

It turns out that the whole analysis of Section 3.4, which ultimately led to
Lemma~\ref{LEM:COND} and Theorem~\ref{THM:GOOD-COND},  can be repeated with obvious 
modifications. Here,
we understand by locally optimal controls $(\ov u,\uzbar)\in \uad^2$ for \CP\ in the sense of 
$L^\infty(0,T)\times\Liq$  pairs for which there is some $\varepsilon>0$ such that
$\,\,\J_2(\S((\jnew{\widehat z}\, \ov u,\uzbar)),(\ov u,\uzbar))\le 
\J_2(\S((\jnew{\widehat z}\,v,v_2)),(v,v_2))\,\,$ for all 
\,$(v,v_2)$\, in \,$\uad^2$ \,with  $\,\|(v,v_2)-(\ov u,\uzbar)\|_{L^\infty(0,T)\times\Liq}\le
\varepsilon.$ We then infer from Theorem~\ref{THM:GOOD-COND} 
\jnew{and the considerations in Remark~\ref{REM2}} the following result.
\Bthm
\label{OP2}
Suppose the conditions of Theorem~\ref{THM:REGU} are satisfied with the exception of \eqref{u:strong}.
Moreover, let \pcol{\ref{b:coeff}}--\ref{G:term}, \eqref{unglatt}, and \eqref{Sce2.1}--\eqref{lumpi1} be
fulfilled. In addition, let $(\ov u,\uzbar)\in\uad^2$ be a locally optimal control for \CP\
in the sense of $L^\infty(0,T)\times \Liq$ \jnew{with associated state $(\ov \mu,
\ov \ph, \ov \s)$ and adjoint state $(p,q,r)$. Then  there are $\ov\lambda_1
\in \partial \G_1(\ov u)$ and $\ov\lambda_2\in \partial \G_2(\uzbar)$ such that the variational inequality \eqref{VUG2} and the pointwise projection formula \eqref{Pro2}
are valid. Moreover, there hold the variational inequality
\begin{align}
\label{VUG3}
&\int_0^T \Bigl(\iO \bigl (-\h(\bvp(x,t))\,\widehat z(x)\,p(x,t) \bigr)\,dx +\kappa\ov\lambda_1(t)+
b_3 \,\ov u(t)\Bigr) \bigl(u(t)-\ov u(t)\bigr)\,dt  \,\ge\, 0 \nonumber\\[1mm]
&\quad\mbox{for all \,$u\in L^\infty(0,T)$\, such that $\,\underline u_1\le u\le 
\widehat u_1$\, a.e. in  $(0,T)$},
\end{align}
as well as for almost every $t\in (0,T)$  the pointwise projection formula}
\begin{align}
\jnew{\ov u(t)=\max\biggl\{\underline u_1(t), \min\biggl\{ \widehat u_1(t),
-b_3^{-1} \biggl(\int_\Omega \bigl(-\h(\bvp(x,t))\widehat z(x) p(x,t) \bigr)dx +\kappa\ov\lambda_1(t) \biggr)\!
\biggr\}\!\biggr\} .}
\end{align}%
\Ethm

\noindent
\underline{{\em Scenario 3:}}
\quad Let again the assumptions of Theorem~\ref{THM:REGU}, with the exception of \eqref{u:strong}, and
\eqref{unglatt} be satisfied. We then consider controls $u_1$ of the general form
\Beq
\label{affine}
u_1=\widetilde z+ {\cal H}[w_1]\,,
\Eeq
with a fixed $\widetilde z\in\Liq\cap L^2(0,T;V)$, and 
where ${\cal H}:L^2(Q) \to\Lzq$  denotes a linear mapping which is continuous from
$\Lzq$ into $\Lzq$ and from $\Liq$ into $\Liq\cap L^2(0,T;V)$. Notice that the 
Scenario 1 studied above is a special case of this one: indeed, in Scenario 1 we had
$\widetilde z=\uebar$ and ${\cal H}=$ null operator. Typical cases for linear operators having
the required regularity properties are spatial convolution integral operators with
sufficiently smooth kernels or the solution operators of a wide class of 
linear parabolic initial-boundary
value problems.

In the above setting, $w_1$ is the true first control variable. We therefore consider as 
set of admissible controls 
\Beq
\label{uadsc3}
\uad^3:=\bigl\{(w_1,u_2)\in \Liq\times\Liq: \ \pcol{\underline u_1\le w_1\le \widehat u_1 \,, \, \ \underline u_2\le u_2\le \widehat u_2 }\, \mbox{ a.e. in }\,Q\bigr\}
\Eeq
under the assumptions \eqref{thresholdreg} on the treshold functions. Now let some $R>0$
be fixed such that 
\Beq
\label{defRsc3}
\uad^3\subset \UR:=\bigl\{(w_1,u_2)\in\Liq\times\Liq: \|(w_1,u_2)\|_{\Liq\times\Liq}<R\bigr\}\,.
\Eeq
Then it follows from the continuity of ${\cal H}$ as a mapping from $\Liq$ into
$\Liq\cap L^2(0,T;V)$ that there exists some $\widehat R>0$ such that for every 
$(w_1,u_2)\in \UR$ it holds 
$\,\,\|(\widetilde z+{\cal H}[w_1],u_2)\|_{(\Liq\cap L^2(0,T;V))\times\Liq}\,<\widehat  R$. 
Therefore, the uniform  separation condition \eqref{separation} is valid for all second solution components associated with controls $(w_1,u_2)\in\UR$, and thus also \eqref{global1}, 
so that the analysis
 carried out in  Section 3 can be repeated analogously. 
 
The cost functional \jnew{is now assumed in} the form\pcol{%
 \[\J_3((\mu,\vp,\s),(w_1,u_2)= J_3((\mu,\vp,\s),(w_1,u_2){)+\kappa \G((w_1,u_2))}\,,\]
with the differentiable part}
\begin{align}
 \label{costduebis}
 J_3((\mu,\vp,\s),(w_1,u_2))=\ &\frac{b_1}2\int_Q|\vp-\widehat\vp_Q|^2+\frac{b_2}2\iO
 |\vp(T)-\widehat \vp_\Omega|^2 \nonumber\\
 &{}+\frac{b_3}2\int_Q(|w_1|^2+|u_2|^2)\,,
 \end{align}
 where it holds $(\mu,\vp,\s)=\S((\widetilde z+{\cal H}[w_1],u_2))$. 
\jnew{Again, we assume the sparsity term in additive form, namely
 \begin{align}
 \label{lumpi2}
 {\cal G}((w_1,u_2))={\cal G}_1(w_1)+{\cal G}_2(u_2)  \quad\mbox{for 
 \,$(w_1,u_2)\in L^2(Q)\times L^2(Q)$,}
 \end{align}
with nonnegative, continuous and convex functionals ${\cal G}_i:L^2(Q)\to\erre$, for
$i=1,2$. Accordingly, we then have}
\begin{align}
\jnew{\partial {\cal G}((w_1,u_2))\,=\,\{(\lambda_1,\lambda_2): \lambda_1\in \partial {\cal G}_1
(w_1) \,\,\,\mbox{and}\,\,\, \lambda_2\in \partial {\cal G}_2(u_2) \}.} \nonumber
\end{align}
Notice that  the affine mapping defined by the identity \eqref{affine} is 
Fr\'echet differentiable between the spaces 
$\Liq$ and $\Liq\cap L^2(0,T;V)$ with the Fr\'echet derivative ${\cal H}$. 
It  follows from Theorem~\ref{THM:DIFF} and the chain rule that the 
control-to-state operator, which in this situation is given by the operator ${\cal T}$,
where 
\Beq
\label{defTsc3}
{\cal T}((w_1,u_2)):=\S(\widetilde z+{\cal H}[w_1],u_2) \quad\mbox{for }\,\,(w_1,u_2)\in \Liq\times\Liq\,,
\Eeq  
is Fr\'echet differentiable from $\Liq\times\Liq$ into \jnew{$\Y$}. For fixed $(\ov w_1,\uzbar)
\in \Liq\times \Liq$, the derivative $D{\cal T}(\ov w_1,\uzbar)$ is given as follows:
for all directions ${\bf h}=(h_1,h_2)\in \Liq\times \Liq$, 
the value $D{\cal T}(\ov w_1,\ov u_2)[(h_1,h_2)]$ is given by the unique solution $(\etah,\psih,\xih)$ to the linearized system 
\eqref{ls1}--\eqref{ls5}, in which in this case the sum of the last two summands in the second line of Eq.~\eqref{ls1}  has to be replaced by the expression 
$$ 
-\iO \h(\bvp)\,{\cal H}[h_1]\,v -\iO \h'(\bvp)(\widetilde  z+{\cal H}[\ov w_1]) \,\psi\,v\,.
$$
Consequently,  the penultimate term in the second line of Eq.~\eqref{adj2} in the adjoint system 
must be replaced by   $\,\,-\h'(\bvp)(\widetilde z+{\cal H}[\ov w_1])\,p$. 
Accordingly, the sum of the
corresponding last two summands in the second line of the calculation \eqref{Rechnung},
which now is performed for ${\bf h}=(h_1,h_2)=(w_1-\ov w_1,u_2-\uzbar)$,
becomes 
$$  
\int_Q \bigl(\h(\bvp)\,{\cal H}[h_1] \,+\,\h'(\bvp)\,(\widetilde z+{\cal H}[\ov w_1])\psi\bigr)p\,.
$$
At this point, it is useful to introduce the dual operator of ${\cal H}$, which is the operator
${\cal H}^*\in {\cal L}(\Lzq,\Lzq)$ defined through the identity
$$
\int_Q  {\cal H}^*[y] \,z=\int_Q y \,{\cal H}[z] \quad\mbox{for all }\,y,z\in\Lzq.
$$  
We then have 
$$
\int_Q \h(\bvp)\,{\cal H}[h_1]\,p = \int_Q {\cal H}^*[\h(\bvp)p] \,h_1\,,
$$
and from the calculation \eqref{Rechnung} \pier{we can} conclude the following result.
\Bthm
\label{OPT3}
Suppose that the assumption of Theorem~\ref{THM:REGU} be satisfied with the exception of
\eqref{u:strong}. 
Moreover, let \pcol{\ref{b:coeff}}--\ref{G:term}, \eqref{unglatt}, and \eqref{affine}--\eqref{lumpi2} be
fulfilled. In addition, assume that $(\ov w_1,\uzbar)\in\uad^2$ is a locally optimal control for \CP\
in the sense of $\Liq\times \Liq$. Then \jnew{there are  $\ov\lambda_1\in
\partial {\cal G}_1(w_1)$ and $\ov \lambda_2\in \partial {\cal G}_2(u_2)$ such that
the variational inequality \eqref{VUG2} and the pointwise projection formula 
\eqref{Pro2} are satisfied. In addition, there  hold the variational inequality 
\begin{align}
\label{VUG4}
&\int_Q \bigl(-{\cal H}^*[\h(\bvp)\,p] +\kappa\ov\lambda_1+
b_3 \,\ov w_1\bigr)(w_1-\ov w_1) \,\ge \,0 \nonumber\\
&\quad\mbox{for all \,$w_1\in\Liq$\, with \,$\underline u_1\le w_1\le \widehat u_1\,$
a.e. in } \,Q,
\end{align}
as well as almost everywhere in $Q$ the pointwise projection formula
\begin{align}
\ov w_1 \,=\, \max \left\{ \underline u_1, \,\min \left\{ \widehat u_1, \,
-b_3^{-1} \bigl(-{\cal H}^*[\h(\ov\ph)p] +\kappa\ov\lambda_1
\bigr)
\right\} \right\}.
\end{align}
}
\Ethm

\section{Some remarks on sparsity}
\label{SPARSITY}
\setcounter{equation}{0}

In this final section, we discuss the concept of sparsity, that is, the possible 
occurrence of proper subregions of the space-time cylinder $Q$ in which 
(locally) optimal controls vanish. The occurrence of sparsity is a consequence of
the choice of the functional $\G$ and the variational inequality \eqref{VUG} through
the form of the associated subdifferential $\partial\G$. In the following, we confine 
ourselves to the concept of {\em full sparsity}, which is connected with the functional
$\G$ given in \eqref{formg} that we consider in the following. Other concepts of
sparsity, like {\em directional sparsity} with respect to time or to space 
(for these concepts see, e.g., \cite{ST}), could also be considered.      

As far as the underlying double-well potential is concerned,  we confine ourselves to the
case of regular potentials, that is, we assume that the assumptions of Theorem~\ref{THM:WEAK},
\pcol{\ref{b:coeff}}--\ref{G:term}, \eqref{glatter}, and \eqref{defUreg}--\eqref{defRreg} are fulfilled. 
We notice, however, that results resembling that of Theorem~\ref{THM:SPARS} below can also be established 
for singular potentials in any of the scenarios considered above in Section 4.

Now
 let $\ubar=(\uebar,\uzbar)
\in \uad$ be a fixed locally optimal control for \CP\ in the sense of $\Liq\times\Liq$
with associated state 
$(\ov\mu,\bvp,\ov\s)=\S(\ubar)$ and adjoint state $(p,q,r)$. Then the statement of
Theorem~\ref{THM:GOOD-COND} is valid. We now introduce the nonnegative, continuous and convex
functional
\Beq
\label{defj}
j:L^2(Q)\to\erre, \quad j(u):=\|u\|_{L^1(Q)} =\int_Q |u|\,.
\Eeq
It is well known (see, \pier{e.g.,} \cite{IT}) that  for any  $u\in\Lzq$ the elements $\lambda$ of the subdifferential $\partial j(u)$ are of the following form: for a.e. $(x,t)\in Q$ it holds
\begin{equation}\label{defdj}
\pier{\lambda(x,t)\ 
\begin{cases}
\ = - 1   & \text{if } \, u(x,t)<0,\\[4pt]
\ \in [-1,+1],  & \text{if } \, u(x,t)=0,\\[4pt]
\ = + 1   & \text{if }\,  u(x,t)> 0.
\end{cases}}
\end{equation}
Next, we observe that we have, for every $\bu=(u_1,u_2)\in L^2(Q)\times L^2(Q)$,
\Beq
\jnew{\G}(\bu)= j(I_1(\bu))+j(I_2(\bu))\,,
\Eeq 
with the projection operators $I_i$, $i=1,2$, that \pier{were} introduced in 
Remark~\ref{REM2}.
As discussed there, the variational inequality \eqref{VUG} then decouples into the two
independent variational inequalities \eqref{VUG1} and \eqref{VUG2} with $\ov\lambda_i
\in \partial j(\ov u_i)$, for $i=1,2$. Consequently, the pointwise projection formulas 
\eqref{Pro1} and \eqref{Pro2} are valid with suitable $\ov\lambda_i\in
\partial j(\ov u_i)$, for $i=1,2$. We then obtain the following (rather standard) sparsity result.
\Bthm
\label{THM:SPARS}
Assume that the \betti{conditions} of Theorem~\ref{THM:WEAK}, \pier{along with}
\pcol{\ref{b:coeff}}--\ref{G:term}, \eqref{glatter}, \eqref{defUreg}--\eqref{defRreg}, are fulfilled.
In addition, suppose that the thresholds functions $\underline u_i, \widehat u_i$, $i=1,2$, 
are constants satisfying the sign condition
\Beq
\label{Vorzeichen}
\underline u_i<0<\widehat u_i, \mbox{\, for \,} i=1.2\,. 
\Eeq
If   $\ubar=(\uebar,\uzbar)
\in \uad$ is a locally optimal control for \CP\ in the sense of $\Liq\times\Liq$
with associated state 
$(\ov\mu,\bvp,\ov\s)=\S(\ubar)$ and adjoint state $(p,q,r)$, then, for almost every $(x,t)\in Q$,
we have the equivalence relations 
\begin{align}
\label{sparsityuno}
&\uebar(x,t)=0  \quad\Longleftrightarrow \quad  |\h(\bvp(x,t))p(x,t)|\,\le\,\kappa\,,\\
\label{sparsitydue}
&\uzbar(x,t)=0 \quad\Longleftrightarrow \quad |r(x,t)|\,\le\,\kappa\,.
\end{align}
\Ethm
\Bdim
We only show \eqref{sparsityuno}; the proof of \eqref{sparsitydue} is analagous. First, we 
have for almost every $(x,t)\in Q$: if $\uebar(x,t)=0$ then, by virtue of the pointwise projection 
formula \eqref{Pro1}, $0=-\h(\bvp(x,t))p(x,t)+\kappa\ov\lambda_1(x,t)$, where
$\ov\lambda_1(x,t) \in [-1,1]$. This obviously implies that \,\,$|\h(\bvp(x,t))p(x,t)|\le\kappa$.

Conversely, we can argue for almost every $(x,t)\in Q$ as follows: suppose that it holds 
$\,\,|\h(\bvp(x,t))p(x,t)|\le\kappa$. If $\uebar(x,t)>0$ \,then\, $\ov\lambda_1(x,t)=1$, and it 
follows from \eqref{Pro1} that \,\,$-b_3^{-1}(\h(\bvp(x.t))p(x,t)+\kappa)>0$. But then 
\,\,$\h(\bvp(x.t)p(x,t)+\kappa<0$, and we obtain that  $\,\,|\h(\bvp(x,t)p(x,t)|
=-\h(\bvp(x,t))p(x,t)>\kappa$, a contradiction. Similar reasoning shows that also the
assumption $\,\uebar(x,t)<0\,$ leads to a contradiction. We therefore must have
$\,\uebar(x,t)=0$.
\Edim

\Brem
\label{REM3}
At this point, it makes sense to ask whether there exists some uniform constant
$\widehat\kappa>0$ such that all locally optimal controls for \CP\ vanish whenever the
sparsity parameter $\kappa$ satisfies $\kappa>\widehat\kappa$. According to the
conditions \eqref{sparsityuno} and \eqref{sparsitydue}, this will be the case if the 
adjoint variables $p$ and $r$ are bounded by a global constant.  And indeed,
this is an immediate consequence of the estimates performed in the proof of Theorem~\ref{THM:ADJ}.
In particular, it was shown there that there is a constant $C>0$, which does not
depend on the special choice of $\ubar\in\UR$, such that (\pier{cf.}~\eqref{adesti7} and
\eqref{adesti2025})
$$
\|r\|_{\Liq} + \|p\|_{\Liq} \,\le\,C.
$$
\Erem


\section*{Acknowledgments}
\betti{PC and ER gratefully mention their affiliation
to the GNAMPA (Gruppo Nazionale per l'Analisi Matematica, 
la Probabilit\`a e le loro Applicazioni) of INdAM (Isti\-tuto 
Nazionale di Alta Matematica) and their collaboration,
as Research Associate, to the IMATI -- C.N.R. Pavia, Italy.
PC and ER also acknowledge the support of Next Generation EU Project No.P2022Z7ZAJ (A unitary mathematical framework for modelling muscular dystrophies).}


\End{document}